\documentclass[12pt]{article}
\usepackage{caption}
\usepackage{color}
\usepackage{eurosym}
\usepackage{lscape}
\usepackage{afterpage}
\usepackage{pstricks}
\usepackage{pst-plot}
\usepackage{longtable}
\usepackage{dcolumn}
\usepackage{pst-node}
\usepackage[dvips]{graphicx}
\usepackage{amsmath}
\usepackage{amssymb}
\usepackage{amsthm}
\usepackage{multirow}
\usepackage{colortab}
\usepackage{color}
\usepackage{array}
\usepackage{slashbox}
\usepackage{colortbl}
\usepackage{shadbox}
\usepackage{textcomp}
\usepackage{pstricks, pst-node}
\usepackage{pst-all}
\usepackage{subfigure}
\usepackage{lscape}
\psset{arrows=->, labelsep=3pt, mnode=circle}

\usepackage{eurosym}

\newfont{\rams}{msbm10 scaled\magstep1}
\newcommand{\rea}{\mbox{\rams \symbol{'122}}}

\newcommand{\rio}{\mathbb{R}}
\begin{document}
\bibliographystyle{plain}
\title{Stochastic Multiobjective Acceptability Analysis for the Choquet integral preference model and the scale construction problem}

\author{ \hspace{0,1cm} Silvia Angilella \thanks{Department of Economics and Business, University of Catania, Corso Italia 55, 95129  Catania, Italy, e-mails: \texttt{angisil\string@unict.it}, \texttt{salvatore.corrente\string@unict.it},  \texttt{salgreco\string@unict.it}},
\hspace{0,1cm} Salvatore Corrente$^*$,
\hspace{0,1cm} Salvatore Greco $^*$}

%\date{29 may 07}
\date{\today}

\maketitle

%%%%%%%%%%%%%%%%%%%%%%%%%%%%%%%%%%%%%%%%%%%%%%%%%%%%%%%%%%%%%%%%
%
%                   ABSTRACT
%                                                                                                                     %
%%%%%%%%%%%%%%%%%%%%%%%%%%%%%%%%%%%%%%%%%%%%%%%%%%%%%%%%%%%%%%%%
\begin{abstract} The Choquet integral is a preference model used in Multiple Criteria Decision Aiding (MCDA) to deal with interactions between criteria. The Stochastic Multiobjective Acceptability Analysis (SMAA) is an MCDA methodology used to take into account imprecision or lack of data in the problem at hand. For example, SMAA permits to compute the frequency that an alternative takes the $k$-th rank in the whole space of the admissible preference parameters, e.g. in case evaluations on the considered criteria are aggregated through the weighted sum model, in the space of weights compatible with the preference information supplied by the Decision Maker (DM). In this paper, we propose to integrate the SMAA methodology with the Choquet integral preference model in order to get robust recommendations taking into account the whole space of preference parameters compatible with the DM's preference information. In case the alternatives are evaluated by all the criteria on a common scale, the preference parameters are given by the capacity expressing the non-additive weights, representing the importance of criteria and their interaction. If the criteria are instead evaluated on different scales, besides the capacity, preference parameters include  the common scale on which the evaluations of criteria  have to be recoded to be compared. Our approach permits to explore the whole space of preference parameters  being capacities and common scales compatible with the DM's preference information.

\medskip


\end{abstract}
%%%%%%%%%%%%%%%%%%%%%%%%%%%%%%%%%%%%%%%%%%%%%%%%%%%%%%%%%%%%%%%%%%%%%%%%%%%%%%%%%%%%%%%%%%%%%%%%%%%%%%%%%%%%%%%%%%%%%%%%%%

%%%%%%%%%%%%%%%%%%%%%%%%%%%%%%%%%%%%%%%%%%%%%%%%%%%%%%%%%%%
\section{Introduction}%%
%%%%%%%%%%%%%%%%%%%%%%%%%%%%%%%%%%%%%%%%%%%%%%%%%%%%%%%%%%%%
\label{introduction}

In Multiple Criteria Decision Aiding (MCDA) (see ~\cite{FigGreEhr} for a survey on MCDA), an alternative $a_k$,  belonging to a finite set
of $l$ alternatives $A=\{a_1, a_2,\ldots a_l\}$,
is evaluated on the basis of a consistent family of $n$ criteria
$G=\{g_1, g_2, \ldots g_n\}$ where, $g_i:A\rightarrow\rea$ is an interval scale of measurement. In this paper, we shall consider also the case in which $g_{i} \colon A \rightarrow {\cal I}$, where 
${\cal I}=\left\{\left[\alpha,\beta\right]\subseteq\rea\right\}$ and $g_{i}(a_k)=\left[\alpha,\beta\right]$ meaning that the evaluation of alternative $a_k$ on criterion $g_i$ can be whichever value inside the interval $\left[\alpha,\beta\right]$. Obviously, if $\alpha=\beta$, then the evaluation of $a_k$ on criterion $g_i$ is given in a precise way.\\
From here on, we will use the terms criterion $g_{i}$ or criterion $i$ interchangeably ($i=1, 2, \ldots, n$) and, without loss of generality, we will assume that all the criteria have to be maximized.

In case the evaluation on criterion $i$ are precise, i.e. $g_{i}(a_k)\in\rea$ for all $k=1,2\,\ldots,l$, we  define a marginal weak preference relation as follows:

\vspace{-0.5truecm}

$$ a_k \,\,\text{is at least as good as}\,\, a_h \,\,\text{with respect to criterion}\,\, i \,\, \Leftrightarrow \,\, g_i(a_k) \geq g_i(a_h).$$

\vspace{-0.2truecm}

To give a recommendation for the decision making problem at hand, evaluations of alternatives on all criteria have to be aggregated. For this reason, three main approaches have been proposed:
\begin{itemize}
	\item the Multi-Attribute Utility Theory (MAUT) ~\cite{Keeney76},
	\item the outranking methods \cite{roy96}, among which the most well known are ELECTRE \cite{Roy_electre} and PROMETHEE \cite{Brans_book,Brans84},
	\item the decision rule approach based on induction of logical ``if..., then...'' decision
rules through Dominance-based Rough Set Approach (DRSA, see \cite{greco2001rough, GMS_DRSA, SGM_2009}). 
\end{itemize}

In MAUT, an overall value function $U \colon \rea^n \rightarrow \rea$  with $U(g_{1}(a_k),\ldots,g_{n}(a_k))=U(a_k)$ is defined such that:

\begin{itemize}

\item $ a_k \,\,\text{is indifferent to}\,\, a_h     \,\, \Leftrightarrow   \,\, U(a_k) = U(a_h),$

\item $ a_k \,\, \text{is preferred to} \,\, a_h      \,\, \Leftrightarrow \,\, U(a_k) > U(a_h),$

\end{itemize}

and the principal aggregation model of value function is the multiple attribute additive utility:

$$U(a_k)=  u_1(g_1(a_k))+ u_2(g_2(a_k))+ \ldots + u_n(g_n(a_k)) \,\,\,\, \mbox{with}\,\,\,\,a_k \in A,$$

\noindent where $u_i$ are non-decreasing marginal value functions  for $i=1,2, \ldots, n.$\\
Outranking methods are instead based on a binary relation $S$ defined on $A$, where $a_kSa_h$ means that ``$a_k$ is at least as good as $a_h$''. A preference, an indifference and for some methods also an incomparability relation can be built on the basis of $S$.
 
A basic assumption of both MAUT and outranking methods is the absence of interaction (synergy or redundancy) of criteria, which is very often    an  unrealistic assumption or a too strong simplification. For example, in evaluating sport cars, a Decision Maker (DM) could consider criteria such as maximum speed, acceleration and price. From one side, maximum speed and acceleration are redundant criteria because, in general, speed cars also have a good acceleration. Therefore, even if these two criteria can be very important for a DM liking sport cars, their comprehensive importance is smaller than the sum of the importance of the two criteria considered separately. From the other side, the two criteria maximum speed and price lead to a synergy effect because a speed car having also a low price is very well appreciated. For such a reason, the comprehensive importance of these two criteria should be greater than the sum of the importance of the two criteria considered separately.

Within MCDA, the interaction of criteria has been dealt using non additive integrals the most well known of which are the Choquet integral \cite{choquet} and the Sugeno integral \cite{sugeno1974theory} (see \cite{Grabisch1996,Grabisch_book_greco,Grabisch2008} for a comprehensive survey on the use of non-additive integrals in MCDA; see also \cite{GL1,GL2,gmg2011,grecorindonebip,grecorindonerob} for some recently proposed extensions of non additive integrals useful in MCDA). Referring to the outranking methods, a notable contribution to deal with interaction of criteria has been given for ELECTRE methods in \cite{FGS2009}. The decision rule approach instead can represent interaction of criteria in a natural way (see \cite{GMS2004}).

One of the greatest drawbacks of the Choquet integral is the elicitation of the preference parameters. To deal with this issue, the DM can provide direct or indirect preference information \cite{angilella2004assessing,marichal2000determination}. The direct preference information consists of providing  all the values of parameters while the indirect preference information (see e.g. \cite{jacquet1982assessing}) is based on giving some preference information related to preference among alternatives from which the parameters values are inferred. \\
Recently, an approach based on the determination of  necessary and possible preference relations within the so-called Robust Ordinal Regression (ROR) \cite{figueira2009building,greco2008ordinal,Greco10} has been presented in \cite{angilella2010non}.

In this paper, we have widely extended the work presented from us \cite{Angilella2012} in which we have applied the Stochastic Multiobjective Acceptability Analysis (SMAA) to explore the whole space of parameters compatible with some preference information provided by the DM related to the importance of criteria and to their interaction (for a survey on SMAA methods see \cite{tervonen_figueira}). In this paper, we have integrated the Choquet integral and the SMAA methods considering also preference information related to pairwise comparison of some reference alternatives. Moreover, we have also considered  alternatives whose evaluations on considered criteria may be given in terms of intervals of possible values. 

Another drawback of the Choquet integral is that all evaluation criteria have to be defined on a common scale in order to compare alternatives on different criteria. This is a strong requirement for this methodology because, for example, in the car evaluation problem cited above, the DM should be able to compare the speed of a car with its acceleration estimating, for example, if the maximum speed of $200$ km/h is as valuable as $35,000$ \euro.

This problem is quite well known in literature (see e.g. \cite{ModaveGrabisch}). One solution has been proposed in \cite{angilella2004assessing} where a common scale was searched for through Monte Carlo simulation. In this paper, we shall show how SMAA permits to  define a common  scale on which one can compare evaluations on different criteria. Moreover, since this common scale, in general, is not unique, we shall describe how SMAA can explore the space of all compatible common scales and the corresponding preference orders that are obtained passing from one common scale to another one. The whole methodology permitting to define a common scale is another original contribution of this paper with respect to \cite{Angilella2012}.\\
The paper is organized as follows. In Section 2, we present the
basic concepts relative to  interaction of criteria and to the Choquet integral. In Section 3, we briefly describe the SMAA methods.
An extension of the SMAA method to the Choquet integral preference model is introduced in Section 4 and illustrated  by several different examples in Section 5. Some conclusions and future directions of research are presented in Section 6.

%%%%%%%%%%%%%%%%%%%%%%%%%%%%%%%%%%%%%%%%%%%%%%%%%
\section{The Choquet integral preference model}%%%%
%%%%%%%%%%%%%%%%%%%%%%%%%%%%%%%%%%%%%%%%%%%%%%%%%
\label{Choquet}

%In this section we introduce the Choquet preference model. First of all, we recall some definitions.

Let $2^{G}$ denote the power set of $G$ (i.e. the set of all subsets of $G$).  A set function $\mu :2^{G}\rightarrow [0,1]$ is  called a capacity (fuzzy measure) on $2^G$ if  the following properties are satisfied:

\begin{description}

\item [{1a)}] $\mu (\emptyset )=0$ and $\mu (G) = 1 $ (boundary
conditions),

\item [{2a)}] $\forall \:S \subseteq T \subseteq G,\,\;\mu (S)\leq \mu (T)$
(monotonicity condition).

\end{description}

\noindent A capacity is called additive if $\mu (S\cup T)=\mu (S)+\mu (T)$, for any $S,T
\subseteq G$ such that $S\cap T=\emptyset $. \\
If a capacity is additive, then it is  uniquely determined by the value assigned to the singletons: $\mu(\{1\}),\mu (\{2\})\ldots ,\mu (\{n\})$.
In such case $\forall \,\,T \subseteq G$, it results:
$$\displaystyle \mu(T)= \sum_{i \in T} \mu(\{i\}).$$

\noindent Whenever the capacity is non additive,  one needs a value $\mu(T)$  for every  subset $T \subseteq G$.  More precisely,
we have to define $2^{|G|}-2$ coefficients,
since the values $\mu(\emptyset)=0$ and $\mu(G)=1$ are already known.

The M\"{o}bius representation  of a capacity $\mu$ (see \cite{Rota}) is 
a function $m:2^{G}\rightarrow \mathbb{R}$ (\cite{Shafer}) defined as follows:

\begin{equation*}
\mu (S)=\sum_{T \subseteq S}m(T).
\end{equation*}

\noindent Note that if $S$ is a singleton, i.e. $S=\{i\}$ with $i=1, 2, \ldots,n$, then $ \mu(\{i\}) =m(\{i\})$. \\
Moreover, if $S$ is a couple (non-ordered pair) of criteria, i.e. $S=\{i,j\}$, then $ \mu(\{i,j\})=m(\{i\}) + m(\{j\}) +m(\{i,j\})$.

The M\"{o}bius representation $m(S)$ can be obtained from $\mu(S)$ as follows:

\begin{equation*}
m(S)=\sum_{T\subseteq S}(-1)^{|S-T|}\mu (T).
\end{equation*}

\noindent In terms of M\"{o}bius representation (see \cite{Chate}), properties \textbf{1a)} and \textbf{2a)} are, respectively, restated  as:

\begin{description}
  \item[1b)] $\displaystyle m(\emptyset )=0,\,\,\,\sum_{T\subseteq G}m(T)=1$,
  \item[2b)] $\displaystyle \forall \: i \in G\:\:\mbox{{and} }\:\: \forall R\subseteq G\setminus\left\{i\right\},\,\,\,\sum_{ T\subseteq R}\:m(T\cup\left\{i\right\})\geq 0.$
\end{description}

In a multicriteria decision problem, the relevance of any criterion $g_i \in G$ can be analyzed   considering $g_i$ either as a singleton or  in combination with all other criteria.
As a result,  the importance of a criterion $g_i\in G$ is evaluated  not only by the capacity $\mu (\{i\})$, but also by considering $\mu (T\cup \{i\})$, $ T\subseteq G\setminus \{i\}$, i.e. when the criterion interacts with all other criteria.

Given $a_k\in A$ and $\mu $ a capacity $\mu$ on $2^G$, the
Choquet integral \cite{choquet} is defined as follows:

\begin{equation}
C_{\mu }(a_k)=\overset{n}{\underset{i=1}{\sum }}\left[
\left(g_{(i)}(a_k)\right) -\left( g_{(i-1)}\left( a_k\right) \right)
\right] \mu \left( N_{i}\right),
\end{equation}

where $_{(\cdot)}$ stands for a permutation of the indices  of criteria such that:

$
g_{(1)}\left( a_k\right)  \leq  g_{(2)}\left( a_k\right)  \leq
...\leq  g_{(n)}\left( a_k\right) ,
$ with $N_{i}=\left\{ (i),\ldots,(n)\right\} $, $i=1,2, \ldots ,n,$ and $g_{(0)}=0$.

In terms of  M\"{o}bius representation \cite{Gilboa:1994}, the Choquet integral may be reformulated as follows:

\begin{equation}
C_{\mu}(a_k)=\underset{T\subseteq G}{\sum }m(T)\underset{i \in
T}{\,\min } \,g_{i}\left( a_k\right).
\end{equation}

Within a Choquet integral preference model, the main difficult task is to infer $2^{|G|}-2$ parameters, giving to them an adequate meaning from a decisional point of view. With the aim of reducing the number of parameters to be elicitated and of avoiding a too strict description of the interactions among criteria, 
in \cite{Grabisch:1997} the concept of fuzzy $k$-additive capacity has been introduced.

\noindent A capacity is called $k$-additive if 
$m(T)=0$ for $T \subseteq G$ such that $|T|>k$ and there exists at least one $T\subseteq G$, with $|T|=k$, such that $m(T)>0$. In particular, a $1$-additive capacity is the standard additive capacity.

Within an MCDA context,  it is easier  and more straightforward to consider $2$-additive capacities only, since in such case the DMs have to  express a preference information on positive and negative interactions between two criteria, neglecting possible interactions among three, four and generally $n$ criteria.
Moreover, by considering $2$-additive measures the computational issue of determining the parameters  is weakened since only  $n+\binom{n}{2}$ coefficients have to be assessed; specifically, in terms of M\"{o}bius representation, a value $ m(\{i\}) $ for every criterion $i$   and a value $ m(\left\{ i,j\right\} )$ for every couple of  criteria $\{i,j\}$. For all these reasons, in the following we shall consider $2$-additive capacities only. However, the methodology we are presenting can be applied to generally non-additive capacities.

The value that a $2$-additive capacity $\mu$ assigns to a set $S\subseteq G$ can be expressed in terms of M\"{o}bius representation as follows:

\begin{equation}
\mu (S)=\underset{i\in S}{\sum }m\left( \left\{ i\right\} \right) +\underset{\left\{ i,j\right\}
\subseteq S}{\sum }m\left( \left\{ i,j\right\} \right) ,\,\,\forall S\subseteq G.
\label{Mobius}
\end{equation}

With regard to 2-additive capacities, properties \textbf{1b)} and \textbf{2b)} have, respectively, the following expressions:

\begin{description}
\item[1c)] $\displaystyle m \left( \emptyset \right) =0,\,\,\,{\sum_{i\in G}
}m\left( \left\{ i\right\} \right) +{\sum_{\left\{ i,j\right\}
\subseteq G} }m\left( \left\{ i,j\right\} \right) =1,$
\item[2c)] $
\left\{
\begin{array}{l}
\displaystyle m\left( \left\{ i\right\} \right) \geq 0,\; \forall i \in G,\\[2mm]
\displaystyle m\left( \left\{ i\right\} \right) +{\sum_{j\in T}}\:m\left( \left\{ i,j\right\} \right) \geq 0, \; \forall i\in G
\:\mbox{and}\:\forall\: T \subseteq G \setminus \left\{i\right\},$ $T\neq\emptyset.
\end{array}
\right.
$
\end{description}

In this case, the Choquet integral of $ a_k \in A$ is given by:

\begin{eqnarray}
C_{\mu}(a_k) &=&\underset{ \{i\}\subseteq G }{\sum }m\left( \left\{
i\right\} \right) \left( g_{i}\left( a_k\right) \right)
+\underset{\left\{ i, j\right\} \subseteq G }{\sum }m\left(
\left\{ i,j\right\} \right) \min  \{g_{i}\left( a_k\right) ,
g_{j}\left( a_k\right) \}.\,\,\,\,\,\,\,\label{Choquet_Mobius}
\end{eqnarray}

Finally, we recall the definitions of the importance and interaction indices for a couple of criteria.

\noindent The importance of criterion $i \in G$ expressed by the Shapley value \cite{Shapley} in case of a $2$-additive capacity can be written as follows:

\begin{equation}
\varphi \left( \left\{ i\right\} \right) =m\left( \left\{ i\right\} \right)
+\underset{j\in G \setminus \left\{ i\right\} }{\sum }\frac{m\left( \left\{ i,j\right\} \right) }{2}.
\label{Shapley_Mobius}
\end{equation}

\noindent The interaction index \cite{Murofushi1993} expressing the sign and the magnitude of the interaction in a couple of criteria $\{i,j \}\subseteq G$, in case of a 2-additive capacity $\mu$, is given by:

\begin{equation}
\varphi \left( \left\{ i,j\right\} \right) =m\left( \left\{ i,j\right\} \right).  \label{Murofushi_Mobius}
\end{equation}

%%%%%%%%%%%%%%%%
\section{SMAA}%%
%%%%%%%%%%%%%%%%
\label{SMAA}%%%%
%%%%%%%%%%%%%%%%
Stochastic Multiobjective Acceptability Analysis (SMAA) methods  ~\cite{Lahdelma,Lahdelma_S2} are a family of MCDA methods aiming to get recommendations on the problem at hand taking into account uncertainty or imprecision on the considered data and preference parameters. %The main idea underlying SMAA is to explore through Monte Carlo  simulation the whole set of criteria weights that give the first position or a certain rank to a specific alternative.
Several SMAA methods have been developed to deal with different MCDA problems:  SMAA-2 has been presented in ~\cite{Lahdelma_S2} for ranking problems, SMAA-O  \cite{Lahdelma2003} has been introduced for multicriteria problems with ordinal criteria and SMAA-TRI \cite{Tervonen2009} for sorting problems. Other two recent contributions relating to SMAA and ROR have been presented in \cite{Kadzinski2013} and \cite{Kadzinski2013A}. For a detailed survey on SMAA methods see \cite{tervonen_figueira}.

Since in this paper we consider ranking problems, we only describe the SMAA-2 method.

In SMAA-2 the most commonly used utility function is the linear one:

$$u(a_{k},w)=\sum_{i=1}^n w_{i}g_{i}(a_k).$$

\noindent In order to take into account imprecision or uncertainty, SMAA-2 considers two probability distributions $f_{W}(w)$ and $f_{\chi}(\xi)$, respectively, on $W$ and $\chi$, where $W=\{(w_1,\ldots,w_n) \in \rio^n \colon w_i \geq 0 \,\,\text{and }\,\, \sum_{i=1}^{n} w_{i}=1\}$ and $\chi$ is the evaluation space.% and $\xi\in\chi$.

First of all, SMAA-2 introduces a ranking function relative to the alternative $a_k$:

$$rank(k,\xi,w)=1+\sum_{h\neq k}\rho\left(u(\xi_{h},w)>u(\xi_{k},w)\right),$$

\noindent where $\rho(false)=0$ and $\rho(true)=1$.

\noindent Then, for each alternative $a_{k}$, for each evaluations of alternatives $\xi\in \chi$ and for each rank $r=1,\ldots,l$, SMAA-2 computes the set of weights of criteria for which alternative $a_k$ assumes rank $r$:

$$W_{k}^{r}(\xi)=\left\{w\in W:rank(k,\xi,w)=r\right\}.$$

SMAA-2 is based on the computation of three indices:

\begin{itemize}
\item The rank acceptability index measures for each alternative $a_{k}$ and for each rank $r$ the variety of different parameters compatible with the DM's preference information giving to $a_k$ the rank $r$ and is given by:

$$b_{k}^{r}=\int_{\xi\in \chi}f_{\chi}(\xi)\int_{w\in W_{k}^{r}(\xi)}f_{W}(w)\;dw\;d\xi.$$

$b_{k}^{r}$ gives the probability that alternative $a_k$ has rank $k$ and it is within the range $[0,1]$.

%The most preferable alternatives are the one with the high acceptability indices for the best ranks.

\item The central weight vector describes the preferences of a typical DM giving to $a_{k}$ the best position and is defined as follows:
$$w_{k}^{c}=\frac{1}{b_{k}^{1}}\int_{\xi \in {\chi}}f_{\chi}(\xi)\int_{w\in W^{1}(\xi)}f_{W}(w)w\;dw \;d\xi;$$

\item The confidence factor is defined as the frequency of an alternative to be the preferred one with the preferences expressed by its central weight vector and is given by:
$$p_{k}^{c}=\int_{{\substack{\xi\in \chi:u(\xi_{k},w_{k}^{c})\geq u(\xi_{h},w_{k}^{c}) \\ \forall h=1,\ldots,l}}}f_{\chi}(\xi)\;d\xi.$$

\end{itemize}

In the following, we shall consider also the frequency that an alternative $a_h$ is preferred to an alternative $a_k$ in the space of the preference parameters (weight vectors in case of SMAA-2), i.e.

$$p(a_h,a_k)=\int_{w\in W} f_W(w) \int_{{\substack{\xi\in \chi:u(\xi_{h},w)\geq u(\xi_{k},w)}}} f_\chi(\xi) d\xi\; dw.$$ 

From a computational point of view, the considered indices are evaluated by the multidimensional integrals approximated by using the Monte Carlo method. \\
After performing  Monte Carlo simulations, the aforementioned indices are evaluated simultaneously in order to help the DM to choose the best alternative of the decision problem under consideration.

%%%%%%%%%%%%%%%%%%%%%%%%%%%%%%%%%%%%%%%%%%%%%%%%%%%%%%%%%%%%%%%%%%%%%%%%%%%%%%%%%%%%
\section{An extension of the SMAA method to the Choquet integral preference model}%%%%
%%%%%%%%%%%%%%%%%%%%%%%%%%%%%%%%%%%%%%%%%%%%%%%%%%%%%%%%%%%%%%%%%%%%%%%%%%%%%%%%%%%%
\label{SMAA-Choquet}%%%%%%%%%%%%%%%%%%%%%%%%%%%%%%%%%%%%%%%%%%%%%%%%%%%%%%%%%%%%%%%%
%%%%%%%%%%%%%%%%%%%%%%%%%%%%%%%%%%%%%%%%%%%%%%%%%%%%%%%%%%%%%%%%%%%%%%%%%%%%%%%%%%%%
In this section, we shall describe how integrate the SMAA methodology within the Choquet integral preference model in three different cases:
\begin{description}
\item[case 1)] the evaluations on the criteria are on a common scale and are given in a precise way, that is $g_{i}(a_k)\in\rea$ for all $i$ and for all $k$,
\item[case 2)] the evaluations on criteria are on a common scale, but  they  can be given in an imprecise way, that is $g_{i}(a_k)\in[\alpha_i^{k},\beta_i^{k}]$ and $\alpha_i^{k}\leq\beta_i^{k}$, for some $i$ and for some $k$,
\item[case 3)] the evaluations on the criteria are on different scales and therefore a common scale has to be constructed (for the sake of simplicity in this case we have supposed that   evaluations of alternatives  on the considered criteria are given in a precise way). 
\end{description}
In Section \ref{Choquet}, we have observed that the use of the Choquet integral in terms of M\"{o}bius representation with a 2-additive capacity involves the knowledge of $n+\binom{n}{2}$ parameters. In order to get these parameters, the DM is therefore asked to provide some preference information in a direct or an indirect way. In the context of the Choquet integral preference model, the direct preference information consists of providing  the  capacity involved in the considered method while the indirect preference information consists in asking the DM  some preference information regarding comparisons of alternatives or importance and interaction of criteria from which eliciting  the capacity. Generally, the indirect preference information requires less cognitive effort from the DM, and for this reason it is widely used in MCDA (see for example \cite{jacquet1982assessing,jacquet2001preference} if the preference model is a value function, \cite{mousseau1998inferring,mousseau2000user} if the preference model is an outranking relation and \cite{angilella2010non} if the preference model is the Choquet integral). Notice that the use of the indirect preference information is intrinsic in the decision rule approach \cite{greco2001rough}.

If the preference model is the Choquet integral, the DM can provide the preference information presented in the following, together with its formulation in terms of linear constraints:

\begin{itemize}
\item Comparisons related to importance and interaction of  criteria:
\begin{itemize}
\item criterion $i$ is at least as important as $j$ ($i\succsim j$): $\varphi(\{g_i\}) \geq\varphi(\{g_j\})$;
\item criterion $i$ is more important than criterion $j$ ($i\succ j$): $\varphi(\{g_i\}) >\varphi(\{g_j\})$;
\item criteria $i$ and $j$ have the same importance ($i\sim j$): $\varphi(\{g_i\}) = \varphi(\{g_j\})$;
\item criteria $i$ and $j$ are synergic: $\varphi(\{g_i,g_j\}) >0$;
\item criteria $i$ and $j$ are redundant: $\varphi(\{g_i,g_j\}) <0$.\\
%where $\varphi(\{g_i\})$ and $\varphi(\{g_i,g_j\})$ have been defined in section \ref{Choquet};
\end{itemize}
\item Comparisons between couples or quadruples of alternatives:
\begin{itemize}
\item alternative $a_k$ is at least as good as alternative $a_h$ ($a_k\succsim a_h$): $C_\mu(a_k)\geq C_\mu(a_h)$;
\item alternative $a_k$ is preferred to alternative $a_h$ ($a_k\succ a_h$): $C_\mu(a_k)>C_\mu(a_h)$;
\item alternative $a_k$ and $a_h$ are indifferent $(a_k\sim a_h)$: $C_\mu(a_k)=C_\mu(a_h)$;
\item alternative $a_k$ is preferred to alternative $a_h$ more than alternative $a_s$ is preferred to alternative $a_t$ $((a_k,a_h)\succ^{*}(a_s,a_t))$: $C_\mu(a_k)-C_\mu(a_h)>C_\mu(a_s)-C_\mu(a_t)$;
\item the difference of preference between $a_k$ and $a_h$ is the same of the difference of preference between $a_s$ and $a_t$ ($(a_k,a_h)\sim^{*}(a_s,a_t)$): $C_\mu(a_k)-C_\mu(a_h)=C_\mu(a_s)-C_\mu(a_t)$.\\
%where $C_{\mu}(a)$ has been defined in section \ref{Choquet}.
\end{itemize}
\end{itemize}

For descriptive purposes, we distinguish three sets of constraints:

\begin{itemize}
\item Monotonicity and boundary constraints 
$$
\left.
\begin{array}{l}
m\left( \{\emptyset\} \right) =0,\,\,\,\underset{g_{i}\in G}{\sum }m\left( \left\{ g_{i}\right\} \right) +\underset{\left\{ g_{i},g_{j}\right\} \subseteq G}{\sum }m\left( \left\{ g_{i},g_{j}\right\} \right) =1,\\[1mm]
 m\left( \left\{ g_{i}\right\} \right) \geq 0,\:\forall g_{i} \in G, \\[1mm]
 m\left( \left\{ g_{i}\right\} \right) +\underset{g_{j}\in T}{\sum
}\:m\left( \left\{g_{i},g_{j}\right\} \right) \geq 0,\forall g_{i}\in G
\:\hbox{and} \:\forall\: T \subseteq G \setminus \left\{ g_{i}\right\}, T\neq\emptyset \\
\end{array}
\right\}(E^{MB})
$$
\item Constraints related to importance and interaction of  criteria:
$$
\left.
\begin{array}{l}
\varphi(\{g_i\}) \geq\varphi(\{g_j\}), \quad\text{if} \,\,i\succsim j,\\[1mm]
\varphi(\{g_i\}) \geq\varphi(\{g_j\})+\varepsilon, \quad\text{if} \,\,i\succ j,\\[1mm]
\varphi(\{g_i\}) = \varphi(\{g_j\}), \quad\text{if} \,\,i\sim j,\\[1mm]
 \varphi(\{g_i,g_j\}) \geq \varepsilon, \quad\text{ if criteria}\,\, i \,\,\text{and} \,\,j \,\,\text{are synergic with\,\,} i,j \in G,\\[1mm]
\varphi(\{g_i,g_j\})\leq-\varepsilon,   \quad \text{ if criteria}\,\, i \,\,\text{and} \,\,j \,\,\text{are redundant with\,\,} i,j \in G,\\[1mm]
\end{array}
\right\}(E^{C})
$$
\item Constraints related to comparisons between alternatives:
$$
\left.
\begin{array}{l}
C_{\mu}(a_k) \geq C_{\mu}(a_h), \quad\text{if} \,\,a_k\succsim a_h,\\[1mm]
C_{\mu}(a_k) \geq C_{\mu}(a_h)+\varepsilon, \quad\text{if} \,\,a_k\succ a_h,\\[1mm]
C_{\mu}(a_k) = C_{\mu}(a_h)     \quad\text{ if}\,\, a_k\sim a_h,\\[1mm]
C_{\mu}(a_k)-C_{\mu}(a_h)\geq C_{\mu}(a_s)-C_{\mu}(a_t)+\varepsilon, \quad\text{if} \,\,(a_k,a_h)\succ^*(a_s,a_t),\\[1mm]
C_{\mu}(a_k)-C_{\mu}(a_h)= C_{\mu}(a_s)-C_{\mu}(a_t), \quad \text{ if}\,\, (a_k,a_h)\sim^*(a_s,a_t),\\[1mm]
\end{array}
\right\}(E^{A})
$$
\end{itemize}
where the strict inequalities have been transformed in weak inequalities in $E^{C}$ and $E^{A}$ by adding an auxiliary variable $\varepsilon.$

We shall call  \textit{compatible model}, a capacity whose M\"{o}bius representation  satisfies the set of constraints $E^{DM}=E^{MB}\cup E^{C}\cup E^{A}$ where $E^{C}$ or $E^{A}$ could be eventually empty if the DM does not provide any information on importance and interaction of criteria, or comparison of alternatives, respectively. \\
 In order to check if there exists at least one compatible model, one has to solve the following linear programming problem:

\begin{equation}\label{compatible}
\begin{array}{l}
\;\;\mbox{max}\; \varepsilon=\varepsilon^{*} \;\;\;s.t.\\[1mm]
\;\; E^{DM}.\\[1mm]
\end{array}
\end{equation}

\noindent If $E^{DM}$ is feasible and $\varepsilon^{*}>0$, then there exists at least one model compatible with the preference information provided by the DM. If E$^{DM}$ is infeasible or $\varepsilon^{\ast}\leq 0$, then one can check which is the minimum set of constraints causing the infeasibility using one of the techniques described in \cite{mousseau2003resolving}.

Each compatible model restores all the information provided by the DM but, generally, gives different recommendations regarding the alternatives that are not provided as example by the DM. For this reason, the choice of only one compatible model could be considered arbitrary and meaningless. 
 The SMAA methodology can then be applied in order to take into account the whole set of  models compatible with the preference information provided by the DM computing the above recalled indices i.e. the rank acceptability indices, central weight vectors and frequency of the preference. Notice that within the Choquet integral preference model we do not work with weights but with capacities expressed in terms of M\"{o}bius representations; for this reason,  the M\"{o}bius representation of the central capacity will be the equivalent of the central weight vector in SMAA. Also ROR \cite{figueira2009building,Greco10}, and Non-Additive ROR (NAROR) in case the evaluations on criteria $g_i$ from $G$ are aggregated through the  Choquet integral   \cite{angilella2010non}, take into account the whole set of compatible models. More precisely within NAROR, and, more generally, within ROR, two preference relations are taken in consideration: for all $a_h,a_k \in A$ 
\begin{itemize}
	\item \textit{necessary preference}, which holds whenever alternative $a_h$ is preferred to alternative $a_k$ for all compatible models, 
	\item \textit{possible preference}, which holds whenever alternative $a_h$ is preferred to alternative $a_k$ for at least one compatible model. 
\end{itemize}

Observe that because the computation of the necessary and possible preference relations involve to solve two linear programming problems for each pair of alternatives in $A$, we could use the SMAA methodology to approximate the two preference relations. In fact: 
\begin{itemize}
\item if $a_h$ is necessarily preferred to $a_k$, then the sum of the frequencies of the preference of $a_{h}$ over $a_{k}$ and  of the indifference between $a_{h}$ and $a_{k}$ is 100\%,
\item if the frequency of the preference of $a_h$ over $a_k$ is greater than zero, then $a_h$ is possibly preferred to $a_k$. 
\end{itemize}
The vice versa of these statements are not true because, from one side, even if the frequency of the preference of $a_h$ over $a_k$ is  100\%, then there could exist one non-sampled compatible model for which $a_k$ is preferred to $a_h$; from the other side, even if $a_h$ is possibly preferred to $a_k$, it is possible that for all sampled compatible models $a_{k}$ is at least as good as $a_h$ and therefore the frequency of the preference of $a_h$ over $a_k$ is 0\%. Observe that the larger is the sample of  compatible models, the better the approximation of the necessary and the possible preference relations obtained with SMAA, such that, in case of an enough large sample of compatible models, one can reasonably accept the approximations of SMAA as results of the ROR approach. Another result of ROR that can be approximated through SMAA is the interval of ranking positions of an alternative given by the extreme ranking analysis \cite{Promethee2010}. Also in this case, the larger the set of  compatible models, the more accurate is the approximation.

Now, we shall describe how the SMAA methodology can be integrated with the Choquet integral in each of the three cases above considered. 
\begin{description} 
\item[case 1)] Since the evaluations on  considered criteria are on a common scale, they are given in a precise way, and the preference information provided by the DM are related only to the importance and interaction of criteria defined by the linear constraints in $E^{DM}$, we apply the Hit-and-Run method \cite{smith1984,Tervonen2012} in order to sample some compatible models. 

The Hit-And-Run sampling starts from the choice of one point (the M\"{o}bius representation of one capacity) inside the polytope $E^{DM}$ delimited by the constraints translating the DM's preference information. At each iteration,   a random direction is sampled from the unit hypersphere that, with the considered position, generates a line. Finally, one point inside the segment whose extremes are the intersection of the line with its boundaries and the current point is sampled. 

At  each iteration of the Hit-and-Run method, one can store the sampled compatible model. Then the Choquet integral of every considered alternative is evaluated with respect to the stored compatible model, such that  at each iteration a complete ranking of alternatives is obtained.
 
In this way, at the end of all the iterations, on the basis of the obtained rankings of the considered alternatives, one can compute the frequency of the preference and of the indifference between two alternatives, the M\"{o}bius representation of the central capacity for each alternative that arrived first at least once and also the rank acceptability index of each alternative with respect to all the possible ranking positions.\\ 
%Remember here that this procedure applied to this particular case has been already described in \cite{Angilella2012}.

\item[case 2)] In this case, the criteria evaluations may be given as intervals of possible values, that is, $g_{i}(a_k)$ can be whichever value inside the interval $\left[\alpha_{i}^{k},\beta_i^{k}\right]$ for some $i$ and some $k$. 

At each iteration, one samples an evaluation matrix $M$ whose element $M_{ki}$, $k=1,\ldots,l$ and $i=1,\ldots,n$ is taken in a random way inside the interval $\left[\alpha_i^{k},\beta_{i}^{k}\right]$. Then, one samples a model compatible with the preference information provided by the DM. Observe that the constraints related to preference information in terms of importance and interaction of criteria, i.e. the constraints  $E^C$, are not dependent from the sampled evaluations. Instead, the constraints related to preference information in terms of pairwise comparisons between alternatives, i.e. the constraints  in $E^{A}$, are  dependent from the sampled evaluations. Thus, we distinguish two hypothesis: $E^{A}=\emptyset$ or $E^{A}\neq\emptyset.$ 

In the first hypothesis ($E^{A}=\emptyset$), the DM is not able or does not want to provide any preference information on comparisons of alternatives. Therefore, the compatible model sampled at each iteration is independent from the sampled evaluation matrix. For this reason, in order to sample models compatible with the preference information provided by the DM, one can use the Hit-and-Run method as done in case 1). The only difference with respect to the first case, is that the Choquet integral of each alternative does not depend only from the compatible model sampled in that iteration, but it depends also from the evaluation matrix  that has been sampled at the beginning of that iteration. 

In the second hypothesis ($E^{A}\neq\emptyset$), the DM is able to provide some preference information on comparisons of alternatives. For this reason, one cannot use the Hit-and-Run method because the set $E^{A}$ depends from the sampled evaluation matrix and therefore at each iteration one has to sample a compatible model from a different set of constraints $E^{DM}$. Due to the change of the constraints in $E^{\text{DM}}$ at each iteration, it is possible that for certain sampled evaluation matrices there is no model compatible with the DM's  preference information. For this reason, at each iteration, before sampling a compatible model, one needs to check if the set of constraints $E^{DM}$ is feasible. In this case, one can sample a compatible model and consequently obtain the complete ranking of the considered alternatives by applying the Choquet integral with the sampled evaluation matrix and the sampled compatible model. 

In both these hypotheses, after all the iterations, one can compute the frequency of the preference or indifference in each couple of alternatives, the M\"{o}bius representation of the central capacity and the rank acceptability index for each alternative. The only difference between the two hypotheses is that, in the second one, to compute the frequency of the preference for each couple of alternatives, one has to consider only those iterations in which,  the set of constraints $E^{DM}$ corresponding to the sampled evaluation matrix is feasible.

\item[case 3)] In the third case, we suppose that the evaluations with respect to considered criteria are on heterogeneous scales. For example, in evaluating a sport car, one can consider criteria such as maximum speed, acceleration, price, comfort etc. and each of them has a different scale. In this case, one can not apply directly the Choquet integral to aggregate the preferences of the DM because, a requisite of the method is that all considered criteria have a common scale. 
\end{description}
In order to deal with this drawback, we propose to construct at each iteration a common scale with a procedure having the following steps for each criterion $g_{i}$:
\begin{itemize}
\item sampling uniformly from the interval $[0,1]$, $l'$ different real numbers $x_{1},\ldots,x_{l'}$ supposing that the different evaluations on $g_i$ are $l',$ with $l' \leq l$,
\item ordering the $l'$ numbers in an increasing way, $x_{i(1)}<\ldots<x_{i(l')}$,
\item assigning $x_{i{(h)}}$ to the alternatives  having the $h$-th evaluation, in an increasing order, on $g_{i}$. 
\end{itemize}

Supposing to deal with the aforementioned car evaluation problem, and looking at the evaluations of the considered cars on criterion acceleration shown in Table \ref{diff_scales}, we proceed as follows:

\begin{itemize}
\item Because the evaluations of the $10$ alternatives on criterion acceleration are all different, we sample 10 different real numbers from the interval $[0,1]$. For example, $x_1=0.81$, $x_{2}=0.90$, $x_{3}=0.12$, $x_{4}=0.91$, $x_{5}=0.63$, $x_{6}=0.09$, $x_{7}=0.27$, $x_{8}=0.54$, $x_{9}=0.95$, $x_{10}=0.96$.
\item We order the 10 numbers in an increasing way:
$ x_{(1)}=0.09< x_{(2)}=0.12< x_{(3)}=0.27< x_{(4)}=0.54< x_{(5)}=0.63< x_{(6)}=0.81< x_{(7)}=0.90< x_{(8)}=0.91< x_{(9)}=0.95< x_{(10)}=0.96$.
\item We assign value $x_{(1)}=0.09$ to PEUGEOT 208 1.6 8V, value $x_{(2)}=0.12$ to FIAT 500 0.9 and so on (see the third column of Table \ref{diff_scales}).
\end{itemize}

\begin{table}[!h]
\begin{footnotesize}
\begin{center}
\caption{Car evaluation with respect to the criterion acceleration and the corresponding scale\label{diff_scales}}
\begin{tabular}[t]{|c|c|c|}
\hline
\textbf{Cars}                     & \textbf{Acceleration} &  \textbf{Scale value} \\ \hline
\textbf{PEUGEOT 208 1.6 8V    }    & 10.9   &  0.09  \\
\hline
\textbf{Citroen C3 }  & 13.5      &  0.90 \\
\hline
\textbf{FIAT 500 0.9 }  & 11   & 0.12 \\
\hline
\textbf{SKODA Fabia 1.2  }        & 14.2 & 0.95  \\
\hline
\textbf{ LANCIA Ypsilon 5p   }        & 11.4 & 0.54 \\
\hline
\textbf{ RENAULT Clio 1.5 dCi 90   }        & 11.3  & 0.27 \\
\hline
\textbf{  SEAT  Ibiza ST 1.2 }        & 14.6   &  0.96  \\
\hline
\textbf{  ALFA ROMEO   MiTo 1.3  }        & 12.9  & 0.81 \\
\hline
\textbf{  TOYOTA Yaris 1.5 }        & 11.8  & 0.63  \\
\hline
\textbf{ VOLKSWAGEN Polo 1.2  }    & 13.9 & 0.91  \\
\hline 
\end{tabular}
\end{center}
\end{footnotesize}
\end{table}

The values $x_{i{(r)}}, i=1,\ldots,n$ and $r=1,\ldots,l'$, become the evaluations of the considered alternatives on the different criteria. In this way, evaluations on all  criteria are expressed on the same common scale and therefore, having the M\"{o}bius representation of a capacity compatible with the DM's preferences, one can compute the Choquet integral of all alternatives.

As described in case 2), the sampling of the scale will influence the sampling of compatible models only if the DM will provide some preference information on the comparison of alternatives. For this reason, one can distinguish, again, between  $E^{A}=\emptyset$ and $E^{A}\neq\emptyset$. \\
If the DM does not provide any preference information in terms of comparisons between alternatives ($E^{A}=\emptyset$), then one can proceed exactly as done in  case of imprecise evaluations, that is using the Hit-and-Run for sampling compatible models and obtaining at each iteration a complete ranking of the considered alternatives by applying the Choquet integral with the sampled common scale and the sampled compatible model. 

If the DM provides some preference information in terms of pairwise comparisons of alternatives ($E^{A}\neq\emptyset$), then the sampling of a model compatible with the DM's  preference information depends from the sampled common scale. For this reason, at each iteration one can sample a common scale and a model compatible with it. In this case, different options are possible:

\begin{itemize}
\item  After all iterations, compute the rank acceptability indices, the M\"{o}bius representation of the central capacity for each alternative and the preference matrix. This is analogous to the case of imprecise evaluations presented  previously  because the sampling of a common scale at the beginning of each iteration is conceptually analogous to the sampling of an evaluation matrix.
\item Sampling a certain number of possible common scales $S_{1},\ldots,S_{iter}$, considering the corresponding feasible sets of constraints $E^{DM}_{1},\ldots,E^{DM}_{iter}$ and denoted by $\varepsilon_{1},\ldots,\varepsilon_{iter}$, the solutions of the linear programming problems

\begin{equation}\label{most discriminant}
\left.
\begin{array}{l}
\;\;\mbox{max}\; \varepsilon \;\;\;s.t.\\[1mm]
\;\; E_{1}^{DM}\\[1mm]
\end{array}
\right\},
\ldots\ldots\ldots\ldots,
\left.
\begin{array}{l}
\;\;\mbox{max}\; \varepsilon \;\;\;s.t.\\[1mm]
\;\; E_{iter}^{DM}\\[1mm]
\end{array}
\right\}
\end{equation}

we can present to the DM the scale $S_{k}$ such that $\varepsilon_{k}=\max\left\{\varepsilon_{1},\ldots,\varepsilon_{iter}\right\}$, that is the most discriminant one admitting a representation of DM's preference information in terms of the Choquet integral preference model.

%\item Sampled a certain number of different scales $S_{1},\ldots,S_{iter}$, considering the corresponding feasible sets of constraints $E^{DM}_{1},\ldots,E^{DM}_{iter}$, and denoted by $\varepsilon_{1},\ldots,\varepsilon_{iter}$, the solutions of the problems (\ref{most discriminant}), for each alternative $a_k$ we can consider the scale $S_{k}$ such that $\varepsilon_{k}=\max\left\{\varepsilon_{1},\ldots,\varepsilon_{iter_k}\right\}$ where $S_1,\ldots,S_{iter_k}$ are the scales among $S_{1},\ldots,S_{iter},$ for which the alternative $a_k$ gets the best position in the complete ranking obtained by applying the Choquet integral.
\end{itemize}

After obtaining the most discriminant common scale, the decision aiding process can continue in one of the following ways:
\begin{description}
	\item[-] asking the capacities directly to the DM,
	\item[-]inducing one capacity compatible with the preference information given by the DM \cite{marichal2000determination},
	\item[-]considering the whole set of capacities compatible with the preference information provided by the DM taking into account necessary and possible preference relations using NAROR \cite{angilella2010non},
	\item[-]applying the SMAA methodology as in case the evaluations of alternatives on considered criteria are given on a common scale from the beginning.

\end{description}

%%%%%%%%%%%%%%%%%%%%%%%%%%%%%%
\section{Some examples}%%%%%%%
%%%%%%%%%%%%%%%%%%%%%%%%%%%%%%
The whole methodology presented in the previous section will be illustrated by several didactic examples.
In the following, where not differently stated, we shall consider uniform probability distributions $f_W$ and $f_\chi$, respectively, on $W$ and $\chi$.

\subsection{An example with all criteria expressed on the same common scale}%%%%
\label{ExampleDidactic}

Let us consider a set of 18 alternatives evaluated on the basis of  4 criteria: $G=\left\{g_1,g_2,g_3,g_4\right\}$.
The criteria are on a $[0,20]$ scale and the scores of every alternative are  reported in Table \ref{FixedEvaluationMatrix}.

\begin{table}[h]
\begin{footnotesize}
\caption{Evaluation matrix}
\label{FixedEvaluationMatrix}
\begin{center}
\resizebox{10cm}{!}{\begin{tabular}{|c|||c|c|c|c|c|c|c|c|c|c|}
\hline
 & \multicolumn{9}{c|}{\textbf{Alternatives}}\\
 \hline
 \textbf{Criteria}       &   $a_1$ &   $a_2$ &   $a_3$ &   $a_4$ &  $a_5$ & $a_6$ & $a_7$ & $a_8$ & $a_9$\\
           \hline
           \hline
${g_1}$ & 15 & 7  & 18 & 9   & 12 & 8  & 14  & 8   & 3  \\
${g_2}$ & 12 & 8  & 8  & 16  & 5  & 3  & 20  & 13  & 17 \\
${g_3}$ & 10 & 14 & 4  & 4   & 14 & 7  & 5   & 15  & 2  \\
${g_4}$ & 7  & 16 & 12 & 16  & 14 & 20 & 10  & 6   & 14 \\
\hline
 \hline
 \textbf{Criteria}       &   $a_{10}$ &   $a_{11}$ &   $a_{12}$ &   $a_{13}$ &  $a_{14}$ & $a_{15}$ & $a_{16}$ & $a_{17}$ & $a_{18}$\\
           \hline
           \hline
${g_1}$         & 4  & 16 & 8  & 17 & 8  & 20 & 12 & 14  & 9  \\
${g_2}$         & 20 & 7  & 11 & 12 & 6  & 7  & 4  & 11  & 13 \\
${g_3}$         & 8  & 14 & 5  & 6  & 7  & 4  & 15 & 12  & 12 \\
${g_4}$         & 9  & 10 & 19 & 8  & 19 & 12 & 13 & 9   & 6  \\
\hline
\end{tabular}}
\end{center}
\end{footnotesize}
\end{table}

Firstly, the DM expresses his/her preference information on the importance and interaction of criteria that can be synthesized as follows:
\begin{itemize}
	\item  $g_{1}$ is more important than  $g_{2}$,
	\item  $g_3$ is more important than $g_{4}$,
	\item  there is a positive interaction between criteria $g_1$ and $g_2$,
	\item  there is a positive interaction between criteria $g_2$ and $g_3$,
	\item  there is a negative interaction between criteria $g_2$ and $g_4$.
\end{itemize}

As explained in Section \ref{SMAA-Choquet}, we use the Hit-and-Run procedure considering a number of iterations equal to $2,000,000$.\\
For each iteration, we sample a model compatible with the preference information provided by the DM, and therefore we compute the Choquet integral of each alternative obtaining a complete ranking. \\
At the end of all iterations, we compute the rank acceptability index $b_{k}^{r}$ for each $k,r=1,\ldots,l$ and the M\"{o}bius representation of the central capacity for each alternative $a_k$ that can get the first rank at least once, as shown, respectively, in Tables \ref{rank_accept_1} and \ref{central_weights_1}. In particular, looking at Table \ref{rank_accept_1}, we observe that alternatives $a_3$, $a_4$, $a_9$, $a_{10}$ and  $a_{18}$ can never be ranked first, $a_{11}$ has reached the first position more than all other alternatives ($b_{11}^{1}=26.22$) and $a_{9}$ is instead the alternative that reached almost always the last position in the obtained rankings $(b_{9}^{18}=93.18)$. 

\begin{table}[h]
\begin{center}
\caption{Rank acceptability indices considering precise evaluations on the considered criteria and preference information only on importance and interaction of criteria}\label{rank_accept_1}

\resizebox{17cm}{!}{\begin{tabular}{|c|c|c|c|c|c|c|c|c|c|c|c|c|c|c|c|c|c|c|}
\hline
\textbf{Alt} & $\mathbf{b_k^{1}}$ & $\mathbf{b_k^{2}}$ & $\mathbf{b_k^{3}}$ & $\mathbf{b_k^{4}}$ & $\mathbf{b_k^{5}}$ & $\mathbf{b_k^{6}}$ & $\mathbf{b_k^{7}}$ & $\mathbf{b_k^{8}}$ & $\mathbf{b_k^{9}}$ & $\mathbf{b_k^{10}}$ & $\mathbf{b_k^{11}}$ & $\mathbf{b_k^{12}}$ & $\mathbf{b_k^{13}}$ & $\mathbf{b_k^{14}}$ & $\mathbf{b_k^{15}}$ & $\mathbf{b_k^{16}}$ & $\mathbf{b_k^{17}}$ & $\mathbf{b_k^{18}}$ \\  
\hline
    $\mathbf{a_1}$    & 1.76  & 9.46  & 16.16 & 10.66 & 13.72 & 12.75 & 10.64 & 7.81  & 4.98  & 3.85  & 2.89  & 2.35  & 1.37  & 1.15  & 0.45  & 0.02  & 0.00  & 0.00 \\
    $\mathbf{a_2}$    &10.88 & 6.03  & 6.19  & 8.15  & 9.32  & 9.00  & 7.11  & 9.02  & 7.34  & 8.75  & 5.75  & 5.86  & 2.87  & 1.94  & 1.38  & 0.41  & 0.00  & 0.00 \\
    $\mathbf{a_3}$    & 0.00  & 4.37  & 3.85  & 3.10  & 3.10  & 5.24  & 6.45  & 5.01  & 6.10  & 6.19  & 6.97  & 9.31  & 6.23  & 8.23  & 9.73  & 10.47 & 5.51  & 0.14 \\
    $\mathbf{a_4}$    & 0.00  & 0.03  & 0.07  & 0.17  & 0.45  & 0.72  & 1.02  & 2.10  & 2.69  & 4.70  & 5.26  & 8.49  & 16.93 & 13.06 & 16.65 & 18.44 & 9.21  & 0.00 \\
    $\mathbf{a_5}$    & 7.60  & 11.15 & 12.89 & 9.43  & 8.21  & 7.68  & 6.83  & 9.69  & 6.97  & 5.97  & 5.69  & 3.02  & 2.62  & 1.79  & 0.44  & 0.03  & 0.00  & 0.00 \\
    $\mathbf{a_6}$    & 1.60  & 1.30  & 1.29  & 1.42  & 1.93  & 2.31  & 2.94  & 3.42  & 3.12  & 5.56  & 4.62  & 5.40  & 6.13  & 8.79  & 9.05  & 20.47 & 19.69 & 0.96 \\
    $\mathbf{a_7}$    & 13.77 & 4.54  & 5.00  & 5.72  & 5.00  & 6.24  & 7.81  & 6.31  & 9.52  & 7.66  & 6.99  & 7.30  & 6.34  & 3.49  & 3.25  & 1.01  & 0.05  & 0.00 \\
    $\mathbf{a_8}$    & 3.31  & 3.11  & 3.94  & 4.78  & 7.53  & 10.01 & 10.04 & 10.69 & 9.42  & 9.88  & 11.04 & 5.03  & 3.48  & 3.21  & 2.83  & 1.65  & 0.03  & 0.00 \\
    $\mathbf{a_9}$    & 0.00  & 0.00  & 0.00  & 0.00  & 0.00  & 0.00  & 0.00  & 0.00  & 0.00  & 0.00  & 0.00  & 0.03  & 0.08  & 0.15  & 0.26  & 1.05  & 5.26  & 93.18 \\
    $\mathbf{a_{10}}$ & 0.00  & 0.00  & 0.00  & 0.00  & 0.00  & 0.02  & 0.07  & 0.19  & 0.78  & 1.32  & 2.66  & 2.93  & 6.32  & 7.41  & 6.84  & 12.78 & 53.79 & 4.88 \\
    $\mathbf{a_{11}}$ & 26.22 & 17.56 & 15.43 & 15.32 & 9.66  & 6.37  & 5.69  & 2.29  & 0.98  & 0.33  & 0.12  & 0.03  & 0.00  & 0.00  & 0.00  & 0.00  & 0.00  & 0.00 \\
    $\mathbf{a_{12}}$ & 1.14  & 1.01  & 1.67  & 1.67  & 1.66  & 2.01  & 2.91  & 4.11  & 5.34  & 5.38  & 9.17  & 10.31 & 11.59 & 21.09 & 15.51 & 5.13  & 0.31  & 0.00 \\
    $\mathbf{a_{13}}$ & 0.59  & 4.56  & 6.53  & 8.77  & 7.27  & 5.96  & 7.18  & 7.05  & 8.10  & 9.06  & 6.20  & 5.31  & 5.05  & 5.78  & 4.36  & 6.20  & 1.53  & 0.51 \\
    $\mathbf{a_{14}}$ & 0.67  & 2.43  & 2.24  & 1.82  & 2.63  & 4.50  & 5.45  & 4.66  & 8.05  & 6.73  & 8.51  & 7.40  & 16.40 & 8.98  & 12.33 & 7.08  & 0.12  & 0.00 \\
    $\mathbf{a_{15}}$ & 12.84 & 5.54  & 3.16  & 3.70  & 5.05  & 5.97  & 4.50  & 5.11  & 5.33  & 5.46  & 7.06  & 5.45  & 5.60  & 6.16  & 7.46  & 7.10  & 4.17  & 0.33 \\
    $\mathbf{a_{16}}$ & 3.02  & 8.29  & 10.75 & 10.93 & 8.96  & 6.73  & 6.79  & 7.22  & 8.80  & 6.58  & 5.17  & 5.26  & 3.16  & 3.38  & 3.40  & 1.47  & 0.08  & 0.00 \\
    $\mathbf{a_{17}}$ & 16.60 & 20.65 & 10.81 & 13.48 & 12.62 & 9.63  & 6.48  & 3.74  & 2.85  & 1.69  & 0.84  & 0.35  & 0.22  & 0.03  & 0.00  & 0.00  & 0.00  & 0.00 \\
    $\mathbf{a_{18}}$ & 0.00  & 0.00  & 0.02  & 0.87  & 2.88  & 4.85  & 8.07  & 11.56 & 9.64  & 10.90 & 11.08 & 16.16 & 5.61  & 5.37  & 6.05  & 6.69  & 0.24  & 0.00 \\
\hline
\end{tabular} }
\end{center}
\end{table}

\vspace{1truecm}
%%%%%%%%%%%%%%%%%%%%%%%%%%%%%%%%%%%%%%%%%%%%%%%%%%%%%%%%%%%%%%%%%%%%%%%%%%%%%%%%%%
%%%%%%%%%%%%%%%%%%%%%%%%%%%%%%%%%%%%%%%%%%%%%%%%%%%%%%%%%%%%%%%%%%%%%%%%%%%%%%%%
\begin{table}[h]
\begin{center}
\caption{M\"{o}bius representation of the central capacity for alternatives being the first at least once in the final ranking obtained by the Choquet integral, considering precise evaluations on considered criteria and preference information only on importance and interaction of criteria}\label{central_weights_1}
\resizebox{14cm}{!}{\begin{tabular}{|c|c|c|c|c|c|c|c|c|c|c|}
\hline
\textbf{Alt/M\"{o}bius} & $\mathbf{m(\left\{1\right\})}$ & $\mathbf{m(\left\{2\right\})}$ & $\mathbf{m(\left\{3\right\})}$ & $\mathbf{m(\left\{4\right\})}$ & $\mathbf{m(\left\{1,2\right\})}$ & $\mathbf{m(\left\{1,3\right\})}$ & $\mathbf{m(\left\{1,4\right\})}$ & $\mathbf{m(\left\{2,3\right\})}$ & $\mathbf{m(\left\{2,4\right\})}$ & $\mathbf{m(\left\{3,4\right\})}$ \\
\hline
  $\mathbf{a_{1}}$ &   0.36  & 0.09  & 0.18  & 0.18  & 0.21  & -0.06 & -0.09 & 0.17  & -0.04 & 0.01 \\
  $\mathbf{a_{2}}$ &   0.25  & 0.11  & 0.37  & 0.31  & 0.08  & -0.09 & -0.07 & 0.08  & -0.07 & 0.04 \\
  $\mathbf{a_{5}}$ &   0.12  & 0.12  & 0.19  & 0.17  & 0.05  & 0.10  & 0.09  & 0.06  & -0.08 & 0.17 \\
  $\mathbf{a_{6}}$ &   0.19  & 0.11  & 0.48  & 0.57  & 0.04  & 0.07  & -0.07 & 0.04  & -0.09 & -0.33 \\
  $\mathbf{a_{7}}$ &  0.28  & 0.20  & 0.22  & 0.19  & 0.18  & -0.06 & 0.02  & 0.06  & -0.09 & -0.01 \\
  $\mathbf{a_{8}}$ & 0.25  & 0.10  & 0.50  & 0.21  & 0.08  & -0.16 & 0.04  & 0.15  & -0.05 & -0.11 \\
  $\mathbf{a_{11}}$ &  0.25  & 0.10  & 0.24  & 0.22  & 0.06  & 0.15  & -0.02 & 0.07  & -0.05 & -0.02 \\
  $\mathbf{a_{12}}$ & 0.20  & 0.11  & 0.45  & 0.49  & 0.12  & -0.01 & -0.04 & 0.06  & -0.05 & -0.34 \\
  $\mathbf{a_{13}}$ & 0.43  & 0.07  & 0.23  & 0.13  & 0.24  & -0.17 & -0.06 & 0.07  & -0.03 & 0.08 \\
  $\mathbf{a_{14}}$ &  0.16  & 0.11  & 0.32  & 0.48  & 0.07  & 0.11  & -0.07 & 0.14  & -0.08 & -0.24 \\
  $\mathbf{a_{15}}$ &   0.49  & 0.11  & 0.36  & 0.22  & 0.07  & -0.23 & -0.01 & 0.09  & -0.06 & -0.04 \\
  $\mathbf{a_{16}}$ &   0.14  & 0.09  & 0.47  & 0.22  & 0.04  & -0.02 & 0.16  & 0.04  & -0.06 & -0.07 \\
  $\mathbf{a_{17}}$ &  0.18  & 0.10  & 0.15  & 0.16  & 0.15  & 0.10  & 0.02  & 0.17  & -0.04 & 0.01 \\
  \hline
\end{tabular}}  
\end{center}
\end{table}
%\end{landscape}

\begin{table}[h]
\begin{center}
\caption{M\"{o}bius representation of the barycenter of compatible capacities considering precise evaluations on considered criteria and preference information only on importance and interaction of criteria}\label{barycent_1}
\resizebox{14cm}{!}{\begin{tabular}{|c|c|c|c|c|c|c|c|c|c|}
\hline
$\mathbf{m(\left\{1\right\})}$ & $\mathbf{m(\left\{2\right\})}$ & $\mathbf{m(\left\{3\right\})}$ & $\mathbf{m(\left\{4\right\})}$ & $\mathbf{m(\left\{1,2\right\})}$ & $\mathbf{m(\left\{1,3\right\})}$ & $\mathbf{m(\left\{1,4\right\})}$ & $\mathbf{m(\left\{2,3\right\})}$ & $\mathbf{m(\left\{2,4\right\})}$ & $\mathbf{m(\left\{3,4\right\})}$ \\ 
\hline
0.26 & 0.12 & 0.27 & 0.22 & 0.098 & 0.008 & -0.001 & 0.09 & -0.06 & -0.07 \\
\hline  
\end{tabular}}  
\end{center}
\end{table}

For each pair of alternatives $(a_h,a_k)$ we compute also the frequency of the preference of $a_h$ over $a_k$ as shown in Table \ref{preference_1}. Looking at the results, alternative $a_{11}$ can be still considered the best because it is preferred to all other alternatives with a frequency of at least the $63.31\%$, while $a_9$ can be considered surely the worst alternative because all alternatives are preferred to it with a frequency at least equal to the $94.77\%.$

Computing  the M\"{o}bius representation of the barycenter of compatible capacities shown in Table \ref{barycent_1}, and applying the Choquet integral, we find the following complete ranking of the considered alternatives:

$$a_{11}\succ a_{17} \succ a_{1} \succ a_{5} \succ a_{16} \succ a_{2} \succ a_{7} \succ a_{8} \succ a_{13} \succ a_{15} \succ a_{18} \succ a_{3} \succ a_{14} \succ a_{12} \succ a_{6} \succ a_{4} \succ a_{10} \succ a_{9}.$$

%%%%%%%%%%%%%%%%%%%%%%%%%%%%%%%%%%%%%%%%%%%%%%%%%%%%%%%%%%%%%%%%%%%%%%%%%%%%%%%%%%
%%%%%%%%%%%%%%%%%%%%%%%%%%%%%%%%%%%%%%%%%%%%%%%%%%%%%%%%%%%%%%%%%%%%%%%%%%%%%%%%
\begin{table}[h]
\begin{center}
\caption{Frequency of the preference between pairs of alternatives considering precise evaluations on considered criteria and preference information only on importance and interaction of  criteria}\label{preference_1}
\resizebox{17cm}{!}{\begin{tabular}{|c|c|c|c|c|c|c|c|c|c|c|c|c|c|c|c|c|c|c|}
\hline
 \textbf{Alt/Alt}   & $\mathbf{a_{1}}$ & $\mathbf{a_{2}}$ & $\mathbf{a_{3}}$ & $\mathbf{a_{4}}$ & $\mathbf{a_{5}}$ & $\mathbf{a_{6}}$ & $\mathbf{a_{7}}$ & $\mathbf{a_{8}}$ & $\mathbf{a_{9}}$ & $\mathbf{a_{10}}$ & $\mathbf{a_{11}}$ & $\mathbf{a_{12}}$ & $\mathbf{a_{13}}$ & $\mathbf{a_{14}}$ & $\mathbf{a_{15}}$ & $\mathbf{a_{16}}$ & $\mathbf{a_{17}}$ & $\mathbf{a_{18}}$ \\
\hline
    $\mathbf{a_{1}}$ & 0.00  & 56.46 & 82.31 & 94.54 & 50.14 & 85.66 & 67.36 & 70.21 & 99.94 & 99.94 & 24.69 & 86.75 & 82.95 & 83.25 & 71.03 & 54.40 & 19.54 & 92.90 \\
    $\mathbf{a_{2}}$ & 43.54 & 0.00  & 68.55 & 90.84 & 42.85 & 91.78 & 54.37 & 59.23 & 99.97 & 98.56 & 24.05 & 86.02 & 60.76 & 86.17 & 61.69 & 50.39 & 30.49 & 75.18 \\
    $\mathbf{a_{3}}$ & 17.69 & 31.45 & 0.00  & 70.87 & 30.27 & 63.25 & 26.46 & 36.13 & 99.33 & 84.09 & 13.18 & 58.31 & 28.54 & 54.25 & 13.88 & 33.60 & 15.02 & 47.75 \\
    $\mathbf{a_{4}}$ & 5.46  & 9.16  & 29.13 & 0.00  & 9.26  & 46.67 & 8.76  & 11.70 & 100.00 & 78.89 & 1.17  & 26.52 & 16.96 & 28.28 & 25.23 & 12.64 & 1.93  & 21.43 \\
    $\mathbf{a_{5}}$ & 49.86 & 57.15 & 69.73 & 90.74 & 0.00  & 91.72 & 59.45 & 66.32 & 99.99 & 99.32 & 18.10 & 85.11 & 63.90 & 86.92 & 62.57 & 76.10 & 33.99 & 82.88 \\
    $\mathbf{a_{6}}$ & 14.34 & 8.22  & 36.75 & 53.33 & 8.28  & 0.00  & 24.38 & 20.54 & 99.03 & 73.12 & 4.41  & 35.32 & 28.96 & 10.39 & 33.08 & 10.31 & 8.08  & 28.92 \\
    $\mathbf{a_{7}}$ & 32.64 & 45.63 & 73.54 & 91.24 & 40.55 & 75.62 & 0.00  & 55.11 & 99.80 & 97.72 & 25.14 & 77.58 & 64.36 & 70.14 & 61.07 & 44.18 & 24.53 & 69.11 \\
    $\mathbf{a_{8}}$ & 29.79 & 40.77 & 63.87 & 88.30 & 33.68 & 79.46 & 44.89 & 0.00  & 100.00 & 98.97 & 12.26 & 76.40 & 52.87 & 73.15 & 55.81 & 37.86 & 14.60 & 78.19 \\
    $\mathbf{a_{9}}$ & 0.06  & 0.03  & 0.67  & 0.00  & 0.01  & 0.97  & 0.20  & 0.00  & 0.00  & 5.23  & 0.00  & 0.00  & 1.19  & 0.13  & 0.70  & 0.03  & 0.00  & 0.08 \\
    $\mathbf{a_{10}}$ & 0.06  & 1.44  & 15.91 & 21.11 & 0.68  & 26.88 & 2.28  & 1.03  & 94.77 & 0.00  & 0.01  & 16.18 & 6.90  & 13.10 & 14.26 & 1.54  & 0.01  & 1.74 \\
    $\mathbf{a_{11}}$ & 75.31 & 75.95 & 86.82 & 98.83 & 81.90 & 95.59 & 74.86 & 87.74 & 100.00 & 99.99 & 0.00  & 94.78 & 81.70 & 94.66 & 77.34 & 85.12 & 63.31 & 98.66 \\
    $\mathbf{a_{12}}$ & 13.25 & 13.98 & 41.69 & 73.48 & 14.89 & 64.68 & 22.42 & 23.60 & 100.00 & 83.82 & 5.22  & 0.00  & 28.98 & 36.20 & 37.16 & 19.12 & 7.47  & 33.20 \\
    $\mathbf{a_{13}}$ & 17.05 & 39.24 & 71.46 & 83.04 & 36.10 & 71.04 & 35.64 & 47.13 & 98.81 & 93.10 & 18.30 & 71.02 & 0.00  & 64.55 & 52.97 & 39.93 & 17.03 & 58.90 \\
    $\mathbf{a_{14}}$ & 16.75 & 13.83 & 45.75 & 71.72 & 13.08 & 89.61 & 29.86 & 26.85 & 99.87 & 86.90 & 5.34  & 63.80 & 35.45 & 0.00  & 41.22 & 17.65 & 8.79  & 37.08 \\
    $\mathbf{a_{15}}$ & 28.97 & 38.31 & 86.12 & 74.77 & 37.43 & 66.92 & 38.93 & 44.19 & 99.30 & 85.74 & 22.66 & 62.84 & 47.03 & 58.78 & 0.00  & 40.69 & 24.17 & 54.33 \\
    $\mathbf{a_{16}}$ & 45.60 & 49.61 & 66.40 & 87.36 & 23.90 & 89.69 & 55.82 & 62.14 & 99.97 & 98.46 & 14.88 & 80.88 & 60.07 & 82.35 & 59.31 & 0.00  & 30.80 & 77.23 \\
    $\mathbf{a_{17}}$ &80.46 & 69.51 & 84.98 & 98.07 & 66.01 & 91.92 & 75.47 & 85.40 & 100.00 & 99.99 & 36.69 & 92.53 & 82.97 & 91.21 & 75.83 & 69.20 & 0.00  & 100.00 \\
    $\mathbf{a_{18}}$ &7.10  & 24.82 & 52.25 & 78.57 & 17.12 & 71.08 & 30.89 & 21.81 & 99.92 & 98.26 & 1.34  & 66.80 & 41.10 & 62.92 & 45.67 & 22.77 & 0.00  & 0.00 \\
\hline
\end{tabular}}  
\end{center}
\end{table}

Let us suppose the DM gives the following further preference information in terms of comparisons of  
alternatives: $a_{16}\succ a_{2}$, $a_{3}\succ a_{14}$ and $a_{13}\succ a_{8}$.

Using the procedure described in Section \ref{SMAA-Choquet}, we obtain the rank acceptability indices, the M\"{o}bius representations of the central capacities and the frequency of the preference between pairs of alternatives shown, respectively, in Tables \ref{rank_accept_2}, \ref{central_weights_2} and \ref{preference_2} in the Appendix.\\
Considering this preference information, we get some new results. Looking at Table \ref{rank_accept_2}, alternative $a_9$ is almost always the worst because it has reached the last position in the rank with a frequency of the $96.89\%$, but, differently from the previous case, $a_{11}$ is no more the alternative obtaining the first rank more than all other alternatives. In fact, $a_{15}$ has been the first with a frequency of the $35.98\%$, while $a_{11}$ has been the first with a frequency of the $30.33\%$. Looking at Table \ref{preference_2}, we can see that in the space of compatible preference models, alternative $a_{11}$ is  preferred to $a_{15}$ with a frequency of the 53.63\%. From Table \ref{central_weights_2} we still get that, in average, the main differences in the M\"{o}bius representation of central capacities of $a_{11}$ and $a_{15}$ are the following:
\begin{itemize}
	\item $m(\{1\})=0.35$ for $a_{11}$ being the first and $m(\{1\})=0.52$ for $a_{15}$ being the first,
	\item $m(\{3\})=0.20$ for $a_{11}$ being the first and $m(\{3\})=0.31$ for $a_{15}$ being the first,
	\item $m(\{1,3\})=0.08$ for $a_{11}$ being the first and $m(\{1,3\})=-0.2$ for $a_{15}$ being the first.
\end{itemize}

Also in this case, we have computed the M\"{o}bius representations of the barycenter of the compatible capacities shown in Table \ref{barycent_2} of the Appendix, and again applying the Choquet integral considering the precise evaluations on the considered criteria, we get the following complete ranking of the considered alternatives:

$$a_{11}\succ a_{15} \succ a_{17} \succ a_{1} \succ a_{13} \succ a_{7} \succ a_{3} \succ a_{5} \succ a_{16} \succ a_{2} \succ a_{8} \succ a_{18} \succ a_{4} \succ a_{12} \succ a_{14} \succ a_{6} \succ a_{10} \succ a_{9}.$$

\noindent We observe that, alternative $a_{15}$  is in the second position, while it was in the middle of the ranking obtained considering the barycenter in case there were not preference information on comparisons between alternatives.

%%%%%%%%%%%%%%%%%%%%%%%%%%%%%%%%%%%%%%%%%%%%%%%%%%%%%%%%%%%%%%%%%%%%%%%%%%%%%%%%%%%
\subsection{Considering imprecision in the evaluations on considered criteria }%%%%
\label{Exampleimprecision}%%%%%%%%%%%%%%%%%%%%%%%%%%%%%%%%%%%%%%%%%%%%%%%%%%%%%%%%%
%%%%%%%%%%%%%%%%%%%%%%%%%%%%%%%%%%%%%%%%%%%%%%%%%%%%%%%%%%%%%%%%%%%%%%%%%%%%%%%%%%%
At this stage, we slightly modify the evaluation matrix shown in Table \ref{FixedEvaluationMatrix}, by introducing  imprecision on the evaluation criteria (see Table \ref{Ev_matrix_imprecision}). %In such case, the criteria  range between a minimum and a maximum value.
In this example, we suppose that the evaluations of considered alternatives on each criterion are integer numbers within an interval (for example, the evaluation of $a_1$ on criterion $g_1$  can be 14, 15 or 16), but this is not a specific requirement for our model i.e., in general, we can  sample values from the whole interval. %looking at Table \ref{Ev_matrix_imprecision}, 
We can consider this as a specific probability distribution $f_\chi(\xi)$ concentrating uniformly the mass only on the integers in the interval of evaluations on considered criteria. 
%%%%%%%%%%%%%%%%%%%%%%%%%%%%%%%%%%%%%%%%%%%%%%%%%%%%%%%%%%%%%%%%%%%%%
%%%%%%%%%%%%%%%%%%%%%%%%%%%%%%%%%%%%%%%%%%%%%%%%%%%%%%%%%%%%%%%%
\begin{table}[h]
\begin{center}
\caption{Imprecise evaluations of alternatives on considered criteria}
\label{Ev_matrix_imprecision}
\resizebox{12cm}{!}{\begin{tabular}{|c|||c|c|c|c|c|c|c|c|c|c|}
\hline
 & \multicolumn{9}{c|}{\textbf{Alternatives}}\\
 \hline
 \textbf{Criteria}       &   $a_1$ &   $a_2$ &   $a_3$ &   $a_4$ &  $a_5$ & $a_6$ & $a_7$ & $a_8$ & $a_9$\\
           \hline
           \hline
${g_1}$ & [14,16] & [6,8] & [17,19] & [8,10]   & [11,13] & [7,9]  & [13,15]  & [7,9]   & 3  \\
${g_2}$ & [11,13] & [7,9]  & [7,9]  & [15,17]  & 5  & 3  & [18,20]  & [12,14]  & [16,18] \\
${g_3}$ & [9,11]      & [13,15] & 4  & [3,5]   & [13,15] & [6,9]  & 5   & [14,15]  & 2  \\
${g_4}$ & [6,9]  & [15,17] & [11,13] & [15,17]  & [13,15] & [18,20] & [9,11]  & 6   & [13,15] \\
\hline
 \hline
 \textbf{Criteria}       &   $a_{10}$ &   $a_{11}$ &   $a_{12}$ &   $a_{13}$ &  $a_{14}$ & $a_{15}$ & $a_{16}$ & $a_{17}$ & $a_{18}$\\
           \hline
           \hline
${g_1}$         & 4  & [15,17] & [7,9]  & [16,18] & [7,9]  & [18,20] & [11,13] & [13,15]  & [8,10]  \\
${g_2}$         & [18,20] & 7  & [10,12] & [11,13] & [6,8]  & [6,9]  & 4  & [10,12]  & [12,14] \\
${g_3}$         & [7,9]  & [13,15] & 5  & [5,7]  & [6,9]  & [3,5]  & [14,16] & [11,13]  & [11,13] \\
${g_4}$         & [8,10]  & [9,11] & [18,20] & 8  & [18,20] & [11,13] & [12,15] & [8,10]   & [5,7]  \\
\hline
\end{tabular}}
\end{center}
\end{table}
%%%%%%%%%%%%%%%%%%%%%%%%%%%%%%%%%%%%%%%%%%%%%%%%%%%%%%%%%%%%%%%%%%%%%%%%%%%%%%

%Within the hypothesis that the criteria are imprecise,  we implement again a hit and run sampling with $2,000,000$ simulations. The sampling on the criteria evaluations is discrete even if it could be also considered continuous.

Firstly, we shall take into account only the preference information in terms of importance and interaction of criteria as presented in Section \ref{ExampleDidactic}. As explained in Section \ref{SMAA-Choquet}, the set of constraints $E^{DM}$ is the same during all the iterations and therefore we can apply the Hit-and-Run method to sample $2,000,000$ compatible models. At the end of all iterations we obtain the rank acceptability index of each alternative for each possible rank, the M\"{o}bius representations of the central capacities and the preference matrix shown, respectively, in Tables \ref{rank_accept_IMP_1}, \ref{central_weights_IMP_1} and \ref{preference_IMP_1} in the Appendix. Comparing results in Tables \ref{rank_accept_1} and \ref{rank_accept_IMP_1}, we observe that the imprecision on the evaluations of alternatives with respect to considered criteria has  increased the uncertainty for the ranking position of the alternatives. In fact, in case of  exact evaluations, five alternatives cannot reach the first position, while considering imprecise interval evaluations, only two alternatives $(a_9$ and $a_{10})$ cannot arrive first in the final ranking. Nevertheless, the best alternative should be still chosen from $a_{11},$ $a_{17},$ and $a_{7}$ because they have the greatest rank acceptabilities for the first position  ($b_{11}^{1}=21.06\%$, $b_{17}^{1}=15.98\%$ and $b_{7}^{1}=11.43\%$) while the worst one is almost surely $a_{9}$ because it has the greatest rank acceptability for the last position ($b_{9}^{18}=90.05\%$). Table \ref{preference_IMP_1} supports the idea that $a_{11}$ is the best alternative because it is preferred to all other alternatives with a frequency at least equal to the $57.60\%$. Table \ref{central_weights_IMP_1} gives the M\"{o}bius representation of the central capacities ranking the considered alternatives in the first position. 

As already done previously, we find the barycenter of the M\"{o}bius representations of the central capacities shown in Table \ref{barycent_3} in the Appendix and representing the average preferences of the DM. Because, we sampled uniformly at each iteration an evaluation matrix, the average evaluation matrix is that one in which each element is the middle point of each interval. Applying the Choquet integral considering the average evaluation matrix and the barycenter shown in Table \ref{barycent_3}, we get the following ranking of alternatives:

$$a_{11}\succ a_{17} \succ a_{1} \succ a_{5} \succ a_{2} \succ a_{16} \succ a_{7} \succ a_{8} \succ a_{13} \succ a_{15} \succ a_{3} \succ a_{18} \succ a_{14} \succ a_{12} \succ a_{4} \succ a_{6} \succ a_{10} \succ a_{9}.$$

Introducing the preference information in terms of comparison between alternatives presented in Section \ref{ExampleDidactic}, we find the rank acceptability indices, the M\"{o}bius representations of the central capacities and the preference matrix shown in Tables \ref{rank_accept_IMP_2}, \ref{central_weights_IMP_2} and \ref{preference_IMP_2} reported in the Appendix. Moreover, considering the barycenter of the M\"{o}bius representations of the central capacities and applying the Choquet integral we get the following ranking of considered alternatives:

$$a_{17}\succ a_{7} \succ a_{11} \succ a_{1} \succ a_{13} \succ a_{15} \succ a_{5} \succ a_{3} \succ a_{16} \succ a_{8} \succ a_{18} \succ a_{2} \succ a_{4} \succ a_{12} \succ a_{14} \succ a_{6} \succ a_{10} \succ a_{9}.$$

\newpage

%%%%%%%%%%%%%%%%%%%%%%%%%%%%%%%%%%%%%%%%%%%%%%%%%%%%%%%%%%%%%%%%%%
\subsection{An example with the  criteria expressed on different scales}%%%%
%%%%%%%%%%%%%%%%%%%%%%%%%%%%%%%%%%%%%%%%%%%%%%%%%%%%%%%%%%%%%%%%%%
\label{Cars_scales}%%%%%%%%%%%%%%%%%%%%%%%%%%%%%%%%%%%%%%%%%%%%%%%
%%%%%%%%%%%%%%%%%%%%%%%%%%%%%%%%%%%%%%%%%%%%%%%%%%%%%%%%%%%%%%%%%%

In this section, we deal with a decision making problem in which the evaluation of alternatives on considered  criteria are expressed  on heterogeneous scales. \\
From the city-car segment market, we select ten cars evaluated on the basis of the following criteria: price (in Euro), acceleration (0/100 km/h in seconds), maximum speed (in km/h) and consumption (in l/100km) (see Table \ref{Ev_Cars}).

\begin{table}[h]\centering
\caption{Evaluation matrix}
\begin{center}\label{Ev_Cars}
\resizebox{12cm}{!}{\begin{tabular}[t]{|c|c|c|c|c|}
\hline
\textbf{Cars}                     & \textbf{Price} & \textbf{Acceleration} & \textbf{Max speed} & \textbf{Consumption}  \\
                                  &  \textbf{Euro}             & \textbf{0/100 km/h}   & \textbf{km/h}   & \textbf{l/100km}  \\ \hline
\textbf{PEUGEOT 208 1.6 8V    }    & 17,800       & 10.9                  & 185            & 3.8  \\
\textbf{   e-HDi 92 CV Stop$\&$Start 3p. Allure}            &               &                       &                 &     \\
\hline
\textbf{Citroen C3 }  & 15,750       & 13.5                  & 163             & 3.8  \\
\textbf{1.4 HDi 70 Seduction    }                         &        &      &               &     \\
\hline
\textbf{FIAT 500 0.9 }                 & 15,050       & 11                  & 173           & 4  \\
\textbf{ TwinAir Turbo Street    }  &        &      &               &     \\
\hline

\textbf{SKODA Fabia 1.2  }        & 15,260      & 14.2                  & 172            & 3.4  \\
\textbf{ TDI CR 75 CV 5p. GreenLine   }  &        &      &               &     \\
\hline
\textbf{ LANCIA Ypsilon 5p   }        & 16,300       & 11.4                  & 183            & 3.8  \\
\textbf{  1.3 MJT 95 CV 5p. S$\&$S Gold  }  &        &      &               &     \\

\hline
\textbf{ RENAULT Clio 1.5 dCi 90   }        & 16,050       & 11.3                  & 176            & 4  \\
\textbf{ CV 3p. Dynamique}                  &        &      &               &     \\

\hline
\textbf{  SEAT  Ibiza ST 1.2 }        & 15,700       & 14.6                  & 173            & 3.4  \\
\textbf{  TDI CR Ecomotive }                  &        &      &               &     \\
\hline
\textbf{  ALFA ROMEO   MiTo 1.3  }        & 17,500       & 12.9                  & 174           & 3.5  \\
\textbf{ JTDm 85 CV S$\&$S Progression }                  &        &      &               &     \\
\hline
\textbf{  TOYOTA Yaris 1.5 }        & 17,800       & 11.8                  & 165            & 3.2  \\
\textbf{  Hybrid 5p. Lounge }                  &        &      &               &     \\
\hline

\textbf{ VOLKSWAGEN Polo 1.2  }        & 17,060       & 13.9                  & 173           & 3.4 \\
\textbf{  TDI 5p. BlueMotion 89g }                  &        &      &               &     \\
\hline
 
\end{tabular}}
\end{center}
\end{table}

%%%%%%%%%%%%%%%%%%%%%%%%%%%%%%%%%%%%%%%%%%%%%%%%%%%%%%%%%%%%%%%%
%%%%%%%%%%%%%%%%%%%%%%%%%%%%%%%%%%%%%%%%%%%%%%%%%%%%%%%%%%%%%%%%

The DM supplies the following preference information in terms of importance and interaction of criteria as well as in terms of comparisons of alternatives:

\begin{itemize}
	\item Preference between alternatives: $a_{5}\succ a_{1}$, $a_{7}\succ a_{6}$, $a_{2}\succ a_{3}$,
	\item Comparisons of importance of criteria: $g_{1}\succ g_{2}$, $g_{4}\succ g_{3}$,
	\item Positive interaction between criteria $g_3$ and $g_4$ and negative interaction between criteria $g_2$ and $g_3$.
\end{itemize}

As explained in Section \ref{SMAA-Choquet}, at each iteration we sample a common scale, and, if the set of constraints $E^{DM}$, that in this case will depend from the sampled scale, is feasible, then we  sample a capacity compatible with these constraints. \\
At the end of all the iterations, we shall get the rank acceptability indices, the M\"{o}bius representations of the central capacities for each alternative and the preference matrix shown respectively in Tables \ref{RAI_Diff_Scales_Cars}, \ref{Central_Weights_Diff_scales_Cars} and \ref{Pref_Diff_Scales_Cars}. \\
Looking at Table \ref{RAI_Diff_Scales_Cars}, we observe that car $a_{7}$ is the most preferred by the DM ($b_{7}^{1}=71.82\%$) followed by $a_{4}$, while $a_{3}$ is most frequently the least preferred car ($b_{3}^{10}=40.68\%$) and $a_1$, $a_2$, $a_3$, $a_6$ and $a_{10}$ can never arrive first. Table \ref{Central_Weights_Diff_scales_Cars} gives the M\"{o}bius representations of the central capacities ranking considered alternatives in the first position at least once, while from Table \ref{Pref_Diff_Scales_Cars}, giving the frequency of the preference between pairs of alternatives, we observe that $a_7$ is preferred to all other alternatives with a frequency at least equal to $71.91\%$.

%Since the criteria are expressed on different scales, at each simulation we randomly generate  a  set of common scales respecting the criteria evaluations' order displayed in Table  \ref{Ev_Cars}.

%With respect to each alternative,  we select  the most representative scale among all the available. 
%The selecting rule considered in the example is to adopt  the one discriminating much more the alternatives' utility evaluations.
%Of course, some other selecting rules can be adopted, for example.... 

%For each alternative, some random weights (M\"{o}bius measures) have been generated. In the example, we consider  $2,000,000$ iterations.  
%
%Then, we outperform the optimization procedure, presented in Section ??? , for the number of iterations in which we have found a set of weights compatible with  the  DM's preference information. 
%
%Finally, we evaluate the rank acceptability indices, the central weights and the percentage of preference between alternatives,   that are, respectively, reported in Tables  

%%%%%%%%%%%%%%%%%%%%%%%%%%%%%%%%%%%%%%%%%%%%%%%%%%%%%%%%%%%%%%%%
%%%%%%%%%%%%%%%%%%%%%%%%%%%%%%%%%%%%%%%%%%%%%%%%%%%%%%%%%%%%%%%%
\begin{table}[h]
  \begin{center}
  \caption{Rank acceptability indices sampling simultaneously compatible capacities and scales}
  
   \resizebox{10cm}{!}{\begin{tabular}{|c|c|c|c|c|c|c|c|c|c|c|}
    \hline
     \textbf{Alt}      & $\mathbf{b_k^1}$ & $\mathbf{b_k^2}$ & $\mathbf{b_k^3}$ & $\mathbf{b_k^4}$ & $\mathbf{b_k^5}$ & $\mathbf{b_k^6}$ & $\mathbf{b_k^7}$ & $\mathbf{b_k^8}$ & $\mathbf{b_k^9}$ & $\mathbf{b_k^{10}}$  \\
    \hline
   $\mathbf{a_1}$    & 0     & 0     & 0 & 0.017 & 0.14 & 1.68 & 17.45 & 25.28 & 30.11 & 25.30 \\
   $\mathbf{a_2}$    & 0     & 0     & 5.32 & 10.41 & 21.22 & 52.43 & 9.74 & 0.80 & 0.04 & 0 \\
   $\mathbf{a_3}$    & 0     & 0     & 0     & 0.18 & 0.54 & 1.41 & 5.25 & 22.96 & 28.94 & 40.68 \\
   $\mathbf{a_4}$    & 28.06 & 71.14 & 0.71 & 0.07 & 0.004 & 0     & 0     & 0     & 0     & 0 \\
   $\mathbf{a_5}$    & 0.04 & 0.189 & 17.15 & 20.88 & 39.46 & 14.50 & 5.28 & 1.73 & 0.74 & 0 \\
   $\mathbf{a_6}$    & 0     & 0     & 0.1 & 0.63 & 1.67 & 6.27 & 36.4 & 28.26 & 19.9 & 6.74 \\
   $\mathbf{a_7}$    & 71.82 & 28.09 & 0.07 & 0.002 & 0     & 0     & 0     & 0     & 0     & 0 \\
   $\mathbf{a_8}$    & 0.058 & 0.252 & 15.62 & 43.49 & 21.7 & 13.12 & 3.899 & 1.57 & 0.259 & 0.013 \\
   $\mathbf{a_9}$    & 0.006 & 0.004 & 0.41 & 1.21 & 4.27 & 7.036 & 20.739 & 19.101 & 19.95 & 27.25 \\
   $\mathbf{a_{10}}$ & 0     & 0.31 & 60.59 & 23.08 & 10.96 & 3.52 & 1.22 & 0.26 & 0.02 & 0 \\
   \hline
    \end{tabular}}
  \label{RAI_Diff_Scales_Cars}%
  \end{center}
\end{table}%
%%%%%%%%%%%%%%%%%%%%%%%%%%%%%%%%%%%%%%%%%%%%%%%%%%%%%%%%%%%%%%%%%%%%%
%%%%%%%%%%%%%%%%%%%%%%%%%%%%%%%%%%%%%%%%%%%%%%%%%%%%%%%%%%%%%%%%
\begin{table}[h]
  \begin{center}
  \caption{M\"{o}bius representations of central capacities sampling simultaneously compatible capacities and scales}
    \resizebox{14cm}{!}{\begin{tabular}{|c|c|c|c|c|c|c|c|c|c|c|}
    \hline
   \textbf{Alt/M\"{o}bius} & $\mathbf{m(\left\{1\right\})}$ & $\mathbf{m(\left\{2\right\})}$ & $\mathbf{m(\left\{3\right\})}$ & $\mathbf{m(\left\{4\right\})}$ & $\mathbf{m(\left\{1,2\right\})}$ & $\mathbf{m(\left\{1,3\right\})}$ & $\mathbf{m(\left\{1,4\right\})}$ & $\mathbf{m(\left\{2,3\right\})}$ & $\mathbf{m(\left\{2,4\right\})}$ & $\mathbf{m(\left\{3,4\right\})}$ \\
   \hline
    %$a_1$ & 0.00  & 0.00  & 0.00  & 0.00  & 0.00  & 0.00  & 0.00  & 0.00  & 0.00  & 0.00 \\
    $\mathbf{a_4}$ & 0.14 & 0.21 & 0.19 & 0.14 & 0.077 & 0.001 & 0.22 & -0.09 & -0.01 & 0.096 \\
    $\mathbf{a_5}$ & 0.17 & 0.22 & 0.21 & 0.17 & 0.002 & 0.001 & 0.19 & -0.109 & 0.001 & 0.109 \\
    $\mathbf{a_7}$ & 0.14 & 0.21 & 0.19 & 0.14 & 0.06 & 0.0006 & 0.23 & -0.1 & -0.007 & 0.099 \\
    $\mathbf{a_8}$ & 0.17 & 0.22 & 0.21 & 0.18 & 0.015 & 0.002 & 0.18 & -0.108 & -0.001 & 0.108 \\
    $\mathbf{a_9}$ & 0.18 & 0.15 & 0.13 & 0.51 & 0.012 & 0.008 & -0.0004 & -0.067 & -0.009 & 0.067 \\
   \hline
   \end{tabular}}
  \label{Central_Weights_Diff_scales_Cars}%
  \end{center}
\end{table}%

%%%%%%%%%%%%%%%%%%%%%%%%%%%%%%%%%%%%%%%%%%%%%%%%%%%%%%%%%%%%%%%%
%%%%%%%%%%%%%%%%%%%%%%%%%%%%%%%%%%%%%%%%%%%%%%%%%%%%%%%%%%%%%%%%

\begin{table}[h]
  \begin{center}
  \caption{Frequency of the preference between pairs of alternatives considering a simulation sampling of random capacities and common scales}
   \resizebox{14cm}{!}{ \begin{tabular}{|c|c|c|c|c|c|c|c|c|c|c|}
    \hline
    \textbf{Alt/Alt}       & $\mathbf{a_1}$ & $\mathbf{a_2}$ & $\mathbf{a_3}$ & $\mathbf{a_4}$ & $\mathbf{a_5}$ & $\mathbf{a_6}$ & $\mathbf{a_7}$ & $\mathbf{a_8}$ & $\mathbf{a_9}$ & $\mathbf{a_{10}}$  \\
    \hline
   $\mathbf{a_1}$  & 0 & 1.64 & 58.82 & 0 & 0     & 31.02 & 0 & 0.619 & 48.18 & 0.31 \\
   $\mathbf{a_2}$  & 98.35 & 0     & 100   & 0     & 25.72 & 95.56 & 0     & 24.59 & 91.2 & 11.08 \\
   $\mathbf{a_3}$  & 41.17 & 0     & 0     & 0     & 2.99 & 16.61 & 0     & 2.38 & 36.3 & 0.67 \\
   $\mathbf{a_4}$  & 99.99 & 100   & 100   & 0     & 99.77 & 100   & 28.08 & 99.67 & 99.98 & 99.65 \\
   $\mathbf{a_5}$  & 100   & 74.27 & 97 & 0.22 & 0     & 97.7 & 0.067 & 39.38 & 90.88 & 23.12 \\
   $\mathbf{a_6}$  & 68.97 & 4.43 & 83.38 & 0     & 2.29 & 0     & 0     & 4.97 & 57.32 & 2.21 \\
   $\mathbf{a_7}$  & 99.99 & 100   & 100   & 71.91 & 99.93 & 100   & 0     & 99.91 & 99.99 & 100 \\
   $\mathbf{a_8}$  & 99.38 & 75.4 & 97.61 & 0.32 & 60.61 & 95.02 & 0.082 & 0     & 96.81 & 23.74 \\
   $\mathbf{a_9}$   & 51.81 & 8.79 & 63.69 & 0.01 & 9.11 & 42.67 & 0.007 & 3.18 & 0     & 0.85 \\
   $\mathbf{a_{10}}$ & 99.68 & 88.91 & 99.32 & 0.34 & 76.87 & 97.78 & 0     & 76.25 & 99.14 & 0 \\
   \hline
    \end{tabular}}
  \label{Pref_Diff_Scales_Cars}%
  \end{center}
\end{table}%
%From Table \ref{Central_Weights_Diff_scales_Cars}, one can observe that alternatives $a_3$, $a_4$ and $a_7$ couldn't be never ranked  first.

%In the final step of the example, we decide to show to the DM the most discriminating scale presented 
% being the most discriminating scale of alternative $a_6$
Since there are different common compatible scales, we propose the DM the most discriminant common scale presented 
in Table \ref{Fixed Scale alternative 6}.
After the DM accepts the common scale, we apply SMAA sampling capacities compatible with the preference information provided by the DM, computing  the rank acceptability indices, the M\"{o}bius representations of the central capacities and the preference matrix displayed, respectively, in Tables \ref{RAI_fixedscales_Cars}, \ref{Central_Weights_Alternative6_Cars} and \ref{pref_fixed_scales_Cars}. 
Applying the Choquet integral with respect to the M\"{o}bius representation of the barycenter of the compatible capacities  shown in Table \ref{barycent_scale}, and considering the most discriminating common scale we get the following ranking of the considered alternatives:

$$a_{7}\succ a_{4} \succ a_{5} \succ a_{8} \succ a_{2} \succ a_{10} \succ a_{3} \succ a_{1} \succ a_{9} \succ a_{6}.$$
  
%Observe that the central common scale obtained through Hit and Run can be used independently by the SMAA methodology, using one among the many methodologies proposed to determine a capacity compatible with the representation of preferences by means of Choquet integral \cite{angilella2004assessing,marichal2000determination}.
%%%%%%%%%%%%%%%%%%%%%%%%%%%%%%%%%%%%%%%%%%%%%%%%%%%%%%%%%%%%%%%%%%%%%%%%%%%%%%%%%%%%%%%%
%%%%%%%%%%%%%%%RISULTATI relativi alla SCALA FISSA PIU RAPPRESENTATIVA(ALTERNATIVA 6)%%%
%%%%%%%%%%%%%%%%%%%%%%%%%%%%%%%%%%%%%%%%%%%%%%%%%%%%%%%%%%%%%%%%%%%%%%%%%%%%%%%%%%%%%%% %%%%%%%%%%%%%%%%%%%%%%%%%%%%%%%%%%%%%%%%%%%%%%%%%%%%%%%%%%%%%%%%%%%%%%%%%%%%%%%%%%%%%%%%%

%%%%%%%SCALA Fissa alternativa 6%%%%%%%%%%%%%%%%%%%%%%%%%

\begin{table}[h]
  \begin{center}
  \caption{Evaluations of alternatives on considered criteria expressed on the most discriminating common scale }
   \resizebox{8cm}{!}{\begin{tabular}{|c|c|c|c|c|}
    \hline
  \textbf{Alt} & \textbf{Price} & \textbf{Acceleration} & \textbf{Max speed} & \textbf{Consumption}  \\
        \hline            
  $\mathbf{a_1}$ & 0.0193 & 0.1135 & 0.6644 & 0.5638 \\
  $\mathbf{a_2}$ & 0.7066 & 0.2816 & 0.0717 & 0.5638 \\
  $\mathbf{a_3}$ & 0.9955 & 0.1187 & 0.3358 & 0.0878 \\
  $\mathbf{a_4}$ & 0.7377 & 0.4395 & 0.09 & 0.6603 \\
  $\mathbf{a_5}$ & 0.4261 & 0.2167 & 0.6397 & 0.5638 \\
  $\mathbf{a_6}$ & 0.4802 & 0.2017 & 0.6284 & 0.0878 \\
  $\mathbf{a_7}$ & 0.7105 & 0.8193 & 0.3358 & 0.6603 \\
  $\mathbf{a_8}$ & 0.2811 & 0.2606 & 0.6164 & 0.6138 \\
  $\mathbf{a_9}$ & 0.0199 & 0.248 & 0.0868 & 0.9567 \\
  $\mathbf{a_{10}}$ & 0.3982 & 0.2835 & 0.3358 & 0.6603 \\
    \hline
    \end{tabular}}
  \label{Fixed Scale alternative 6}%
  \end{center}   
 \end{table}%
%%%%%%%%%%%%%%%%%%%%%%%%%%%%%%%%%%%%%%%%%%%%%%%%%%%%%%%%%%%%%%%%%%%%%
%%%%%%%%%%%%%%%%%%%%%%%%%%%%%%%%%%%%%%%%%%%%%%%%%%%%%%%%%%%%%%%%

\begin{table}[htbp]
  \begin{center}
  \caption{Rank acceptability indices taking into account evaluations of alternatives on considered criteria expressed on the most discriminating common scale shown in Table  \ref{Fixed Scale alternative 6}}
    \resizebox{10cm}{!}{ \begin{tabular}{|c|c|c|c|c|c|c|c|c|c|c|}
    \hline
     \textbf{Alt}       & $\mathbf{b_k^1}$ & $\mathbf{b_k^2}$ & $\mathbf{b_k^3}$ & $\mathbf{b_k^4}$ & $\mathbf{b_k^5}$ & $\mathbf{b_k^6}$ & $\mathbf{b_k^7}$ & $\mathbf{b_k^8}$ & $\mathbf{b_k^9}$ & $\mathbf{b_k^{10}}$  \\
    \hline
 $\mathbf{a_1}$  &  0  & 0  & 0.21  & 4.09  & 4.61 & 6.82 & 16.89 & 24.35 & 28.87 & 14.11 \\
 $\mathbf{a_2}$  &  0  & 0  & 27.23 & 13.39 & 13.17 & 21.76 & 15.13 & 8.53 & 0.75 & 0 \\
 $\mathbf{a_3}$  &  0  & 0  & 0     & 7.55 & 5.34 & 6.41 & 20.99 & 19.99 & 25.36 & 14.32 \\
 $\mathbf{a_4}$  &  0  & 65.12 & 10.82 & 10.44 & 6.72 & 5.17 & 1.69 & 0.001 & 0     & 0 \\
 $\mathbf{a_5}$  &  0.52 & 21.71 & 34.57 & 23.30 & 12.25 & 4.92 & 2.17 & 0.53 & 0   & 0 \\
 $\mathbf{a_6}$  &  0    & 0     & 0     & 0     & 0.001 & 0.502 & 7.71 & 28.22 & 27.76 & 35.79 \\
 $\mathbf{a_7}$  &  97.11 & 1.809 & 0.85 & 0.21 & 0.003 & 0     & 0     & 0     & 0     & 0 \\
 $\mathbf{a_8}$  &  1.12 & 7.79 & 19.78 & 21.83 & 17.31 & 20.66 & 10.45 & 1.02  & 0     & 0 \\
 $\mathbf{a_9}$  &  1.22 & 3.43 & 3.17 & 1.98 & 2.7 & 5.57 & 13.58 & 15.42 & 17.12 & 35.76 \\
 $\mathbf{a_{10}}$  &  0     & 0.11 & 3.32 & 17.16 & 37.86 & 28.16 & 11.35 & 1.89 & 0.109 & 0 \\
  \hline
    \end{tabular}}
  \label{RAI_fixedscales_Cars}%
  \end{center}
\end{table}%

%%%%%%%%%%%%%%%%%%%%%%%%%%%%%%%%%%%%%%%%%%%%%%%%%%%%%%%%%%%%%%%%%%%%%
%%%%%%%%%%%%%%%%%%%%%%%%%%%%%%%%%%%%%%%%%%%%%%%%%%%%%%%%%%%%%%%%

\begin{table}[htbp]
  \begin{center}
  \caption{M\"{o}bius representation of the central capacities,  taking into account evaluations of alternatives on considered criteria expressed on the most discriminating common scale, shown in Table \ref{Fixed Scale alternative 6}}\label{Central_Weights_Alternative6_Cars}
  \resizebox{14cm}{!}{\begin{tabular}{|c|c|c|c|c|c|c|c|c|c|c|}
    \hline
  \textbf{ Alt/M\"{o}bius } & $\mathbf{m(\left\{1\right\})}$ & $\mathbf{m(\left\{2\right\})}$ & $\mathbf{m(\left\{3\right\})}$ & $\mathbf{m(\left\{4\right\})}$ & $\mathbf{m(\left\{1,2\right\})}$ & $\mathbf{m(\left\{1,3\right\})}$ & $\mathbf{m(\left\{1,4\right\})}$ & $\mathbf{m(\left\{2,3\right\})}$ & $\mathbf{m(\left\{2,4\right\})}$ & $\mathbf{m(\left\{3,4\right\})}$ \\
   \hline
$\mathbf{a_5}$  &  0.089 & 0.106 & 0.155 & 0.2 & -0.016 & 0.087 & 0.09 & -0.03 & -0.015 & 0.336 \\

$\mathbf{a_7}$  &  0.222 & 0.201 & 0.163 & 0.244 & 0.038 & 0.002 & 0.102 & -0.06 & -0.007 & 0.09 \\

$\mathbf{a_8}$  &  0.145 & 0.098 & 0.167 & 0.28 & 0.017 & -0.017 & -0.027 & -0.033 & -0.012 & 0.38 \\

$\mathbf{a_9}$  &  0.269 & 0.152 & 0.099 & 0.685 & -0.005 & 0.016 & -0.187 & -0.03 & -0.05 & 0.055 \\
\hline
     \end{tabular}}
  \end{center}%
\end{table}%

%%%%%%%%%%%%%%%%%%%%%%%%%%%%%%%%%%%%%%%%%%%%%%%%%%%%%%%%%%%%%%%%
%%%%%%%%%%%%%%%%%%%%%%%%%%%%%%%%%%%%%%%%%%%%%%%%%%%%%%%%%%%%%%%%

\begin{table}[h]
\begin{center}
\caption{ M\"{o}bius representation of the barycenter of capacities taking into account evaluations of alternatives on the considered criteria expressed on the most discriminant common scale shown in Table \ref{Fixed Scale alternative 6}}\label{barycent_scale}
\resizebox{14cm}{!}{\begin{tabular}{|c|c|c|c|c|c|c|c|c|c|}
\hline
$\mathbf{m(\left\{1\right\})}$ & $\mathbf{m(\left\{2\right\})}$ & $\mathbf{m(\left\{3\right\})}$ & $\mathbf{m(\left\{4\right\})}$ & $\mathbf{m(\left\{1,2\right\})}$ & $\mathbf{m(\left\{1,3\right\})}$ & $\mathbf{m(\left\{1,4\right\})}$ & $\mathbf{m(\left\{2,3\right\})}$ & $\mathbf{m(\left\{2,4\right\})}$ & $\mathbf{m(\left\{3,4\right\})}$ \\ 
\hline
0.22 & 0.19 & 0.16 & 0.25 & 0.037 & 0.003 & 0.097 & -0.05 & -0.008 & 0.09 \\
\hline  
\end{tabular}}  
\end{center}
\end{table}

%%%%%%%%%%%%%%%%%%%%%%%%%%%%%%%%%%%%%%%%%%%%%%%%%%%%%%%%%%%%%%%%%%%%%
%%%%%%%%%%%%%%%%%%%%%%%%%%%%%%%%%%%%%%%%%%%%%%%%%%%%%%%%%%%%%%%%
\begin{table}[h]
  \begin{center}
   \caption{Frequency of the preference between pairs of alternatives taking into account evaluations of alternatives on the most discriminant common scale shown in Table \ref{Fixed Scale alternative 6}}\label{pref_fixed_scales_Cars}%
    \resizebox{14cm}{!}{\begin{tabular}{|c|c|c|c|c|c|c|c|c|c|c|}
    \hline
    \textbf{Alt/Alt}       & $\mathbf{a_1}$ & $\mathbf{a_2}$ & $\mathbf{a_3}$ & $\mathbf{a_4}$ & $\mathbf{a_5}$ & $\mathbf{a_6}$ & $\mathbf{a_7}$ & $\mathbf{a_8}$ & $\mathbf{a_9}$ & $\mathbf{a_{10}}$  \\
    \hline
    $\mathbf{a_1}$     &  0   & 18.33 & 48.3 & 6.52 & 0     & 64.7 & 0.21 & 0     & 56.69 & 9.97 \\
    $\mathbf{a_2}$     &  81.66 & 0     & 100   & 0     & 30.88 & 98.49 & 0     & 43.23 & 80.83 & 52.06 \\
    $\mathbf{a_3}$     &  51.69 & 0     & 0     & 0     & 8.66  & 74.71 & 0     & 15.11 & 58.21 & 17.68 \\
    $\mathbf{a_4}$     &  93.47 & 100   & 100   & 0     & 69.57 & 99.99 & 0     & 75.77 & 92.34 & 87.75 \\
    $\mathbf{a_5}$     &  100   & 69.11 & 91.33 & 30.42 & 0     & 100   & 1.27 & 77.77 & 90.54 & 88.35 \\
    $\mathbf{a_6}$     &  35.29 & 1.5 & 25.28 & 0.004 & 0     & 0     & 0     & 0.18 & 47.07 & 0.022 \\
    $\mathbf{a_7}$     &  99.7843 & 100   & 100   & 100   & 98.72 & 100   & 0     & 98.55 & 98.75 & 100 \\
    $\mathbf{a_8}$     &  100   & 56.76 & 84.88 & 24.22 & 22.22 & 99.81 & 1.44 & 0     & 90.09 & 65.22 \\
    $\mathbf{a_9}$     &  43.308 & 19.16 & 41.78 & 7.65 & 9.45 & 52.92 & 1.24 & 9.904 & 0     & 11.73 \\
    $\mathbf{a_{10}}$  &  90.02 & 47.93 & 82.31 & 12.24 & 11.64 & 99.97 & 0     & 34.77 & 88.26 & 0 \\
    \hline
     \end{tabular}}
  \end{center}
\end{table}%

\newpage

\vspace{8truecm}
%%%%%%%%%%%%%%%%%%%%%%%%%
%%%%%%%%%%%%%%%%%%%%%%%%%
\section{Conclusions}%%%%
%%%%%%%%%%%%%%%%%%%%%%%%%
In this paper, we  have integrated the Stochastic Multiobjective Acceptability Analysis (SMAA)  to the Choquet integral preference model extending a work already published by the authors \cite{Angilella2012}. We have proposed to explore the space of the parameters compatible with some preference information provided by the DM using SMAA. In particular, we have considered the DM's preference information not only in terms of relative importance of criteria and interaction between them , but differently from \cite{Angilella2012}, also in terms of pairwise comparison between alternatives and comparisons of intensity of preferences between pairs of alternatives. Moreover, again differently from \cite{Angilella2012}, we have considered also imprecise evaluations of alternatives on the considered criteria expressed in terms of intervals of possible values.

Finally, we have proposed a methodology to construct the common scale required by the Choquet integral; this is very useful in case the criteria for the decision problem at hand are defined on different scales. Such aspect of the methodology we are proposing, constitutes another original contribution with respect to \cite{Angilella2012}. We have provided several different examples in which the proposed methodology has been applied. We envisage the following future developments:
\begin{itemize}
	\item application of SMAA methodology to some extensions of the classical Choquet integral, e.g. the bipolar Choquet integral \cite{GL1,GL2} (see also \cite{GrecoMatarazzoSlowinski02}), the level dependent Choquet integral \cite{gmg2011}, the robust Choquet integral \cite{grecorindonerob};
	\item implementation of the SMAA methodology to the Choquet integral in presence of hierarchy of criteria \cite{acg2013} within the so called multiple criteria hierarchy process \cite{CGShierarchy,CGSoutr}.  
\end{itemize}
We believe that  the methodology we are proposing can greatly contribute to extend and to improve the use of the Choquet integral preference model in Multiple Criteria Decision Aiding.

\bibliographystyle{plain}
\bibliography{General}
%%%%%%%%%%%%%%%%%%%%%%%%%%%%%%%%%%%%%APPENDIX%%%%%%%%%%%%%%%%%%%%%%%%%%%%%%%%%%%%%%%%%%%%%%%%%%%%%%%%%%%%%%%
%\newpage
%\begin{small}
%\subsection*{Appendix}

\newpage
%%%%%%%%%%%%%%%%%%%%%%%%%%%%%%%%%%%%%%
\section*{Appendix}\label{appendix}%%%
%%%%%%%%%%%%%%%%%%%%%%%%%%%%%%%%%%%%%%

%%%%%%%%%%%%%%%%%%%%%%%%%%%%%%%%%%%%%%%%%%%%%%%%%%%%%%
% FIXED EVALUATIONS AND BOTH TYPES OF PREFERENCE%%%%%%
%%%%%%%%%%%%%%%%%%%%%%%%%%%%%%%%%%%%%%%%%%%%%%%%%%%%%%

\begin{table}[h]
\begin{center}
\caption{Rank acceptability indices  considering preference information both in terms of comparison of alternatives and importance and interaction of criteria}\label{rank_accept_2}
\resizebox{17cm}{!}{
\begin{tabular}{|c|c|c|c|c|c|c|c|c|c|c|c|c|c|c|c|c|c|c|}
\hline
\textbf{Alt} & $\mathbf{b_k^{1}}$ & $\mathbf{b_k^{2}}$ & $\mathbf{b_k^{3}}$ & $\mathbf{b_k^{4}}$ & $\mathbf{b_k^{5}}$ & $\mathbf{b_k^{6}}$ & $\mathbf{b_k^{7}}$ & $\mathbf{b_k^{8}}$ & $\mathbf{b_k^{9}}$ & $\mathbf{b_k^{10}}$ & $\mathbf{b_k^{11}}$ & $\mathbf{b_k^{12}}$ & $\mathbf{b_k^{13}}$ & $\mathbf{b_k^{14}}$ & $\mathbf{b_k^{15}}$ & $\mathbf{b_k^{16}}$ & $\mathbf{b_k^{17}}$ & $\mathbf{b_k^{18}}$ \\  
\hline
$\mathbf{a_{1}}$  &    1.25  & 7.05  & 20.65 & 15.61 & 22.70 & 14.64 & 6.70  & 4.49  & 4.33  & 2.02  & 0.49  & 0.04  & 0.02  & 0.01  & 0.00  & 0.00  & 0.00  & 0.00 \\
$\mathbf{a_{2}}$  &   0.00  & 0.02  & 0.89  & 1.72  & 2.86  & 3.63  & 3.53  & 4.40  & 5.90  & 21.34 & 15.26 & 17.32 & 8.53  & 7.24  & 5.72  & 1.61  & 0.03  & 0.00 \\
$\mathbf{a_{3}}$  &   0.00  & 14.80 & 11.54 & 7.51  & 5.84  & 9.47  & 14.00 & 8.32  & 9.18  & 7.55  & 5.30  & 4.85  & 1.00  & 0.50  & 0.14  & 0.00  & 0.00  & 0.00 \\
$\mathbf{a_{4}}$  &    0.00  & 0.00  & 0.00  & 0.03  & 0.10  & 0.28  & 0.42  & 2.95  & 2.10  & 3.77  & 4.97  & 9.02  & 28.28 & 11.80 & 14.35 & 19.25 & 2.67  & 0.00 \\
$\mathbf{a_{5}}$  &    4.06  & 5.81  & 7.42  & 8.10  & 7.92  & 7.97  & 8.70  & 25.07 & 11.06 & 8.34  & 3.53  & 1.52  & 0.46  & 0.03  & 0.00  & 0.00  & 0.00  & 0.00 \\
$\mathbf{a_{6}}$  &    0.00  & 0.00  & 0.01  & 0.03  & 0.02  & 0.04  & 0.12  & 0.24  & 0.42  & 1.42  & 2.63  & 4.17  & 4.58  & 11.66 & 8.52  & 40.58 & 23.45 & 2.12 \\
$\mathbf{a_{7}}$  &    17.21 & 7.07  & 6.83  & 9.79  & 8.15  & 9.44  & 13.66 & 6.32  & 9.19  & 6.23  & 3.85  & 1.32  & 0.62  & 0.28  & 0.03  & 0.00  & 0.00  & 0.00 \\
$\mathbf{a_{8}}$  &    0.00  & 0.00  & 0.01  & 0.07  & 0.54  & 2.18  & 3.95  & 11.22 & 10.55 & 20.45 & 26.58 & 8.63  & 4.71  & 4.69  & 4.12  & 2.31  & 0.01  & 0.00 \\
$\mathbf{a_{9}}$  &    0.00  & 0.00  & 0.00  & 0.00  & 0.00  & 0.00  & 0.00  & 0.00  & 0.00  & 0.00  & 0.00  & 0.00  & 0.01  & 0.01  & 0.04  & 0.27  & 2.78  & 96.89 \\
$\mathbf{a_{10}}$  &    0.00  & 0.00  & 0.00  & 0.00  & 0.00  & 0.02  & 0.06  & 0.13  & 0.37  & 0.50  & 1.09  & 1.00  & 3.18  & 5.37  & 5.19  & 11.73 & 70.38 & 0.99 \\
$\mathbf{a_{11}}$  &    30.33 & 18.57 & 14.55 & 14.01 & 9.74  & 6.68  & 5.74  & 0.34  & 0.04  & 0.01  & 0.00  & 0.00  & 0.00  & 0.00  & 0.00  & 0.00  & 0.00  & 0.00 \\
$\mathbf{a_{12}}$  &    0.00  & 0.00  & 0.02  & 0.05  & 0.10  & 0.18  & 0.27  & 0.66  & 1.08  & 2.05  & 3.82  & 6.98  & 12.93 & 37.34 & 28.29 & 5.88  & 0.36  & 0.00 \\
$\mathbf{a_{13}}$  &    0.98  & 8.39  & 15.35 & 17.00 & 16.57 & 11.90 & 12.89 & 10.13 & 4.71  & 1.24  & 0.39  & 0.30  & 0.12  & 0.02  & 0.00  & 0.00  & 0.00  & 0.00 \\
$\mathbf{a_{14}}$  &    0.00  & 0.00  & 0.00  & 0.00  & 0.04  & 0.08  & 0.11  & 0.37  & 0.91  & 2.85  & 6.24  & 8.11  & 26.64 & 14.16 & 27.83 & 12.35 & 0.31  & 0.00 \\
$\mathbf{a_{15}}$  &    35.98 & 12.28 & 4.75  & 6.14  & 7.22  & 10.41 & 5.66  & 5.24  & 4.32  & 2.77  & 2.75  & 1.16  & 0.58  & 0.34  & 0.33  & 0.06  & 0.00  & 0.00 \\
$\mathbf{a_{16}}$  &   0.45  & 3.89  & 6.58  & 7.90  & 8.06  & 6.90  & 7.72  & 12.19 & 25.64 & 6.21  & 7.45  & 3.33  & 2.49  & 1.13  & 0.06  & 0.00  & 0.00  & 0.00 \\
$\mathbf{a_{17}}$  &   9.73  & 22.11 & 11.40 & 12.01 & 9.81  & 14.51 & 13.98 & 3.13  & 1.97  & 0.87  & 0.45  & 0.02  & 0.01  & 0.00  & 0.00  & 0.00  & 0.00  & 0.00 \\
$\mathbf{a_{18}}$  &  0.00  & 0.00  & 0.00  & 0.04  & 0.33  & 1.69  & 2.50  & 4.81  & 8.22  & 12.38 & 15.18 & 32.23 & 5.83  & 5.42  & 5.39  & 5.98  & 0.01  & 0.00 \\
\hline
\end{tabular}}  
\end{center}
\end{table}

%%%%%%%%%%%%%%%%%%%%%%%%%%%%%%%%%%%%%%%%%%%%%%%%%%%%%%%%%%%%%%%%%%%%%%%%%%%%%%%%
%%%%%%%%%%%%%%%%%%%%%%%%%%%%%%%%%%%%%%%%%%%%%%%%%%%%%%%%%%%%%%%%%%%%%%%%%%%%%%%%
\begin{table}[h]
\begin{center}
\caption{M\"{o}bius representations of central capacities for alternatives taking into account precise evaluations on considered criteria and preference information both in terms of comparison of alternatives and importance and interaction of criteria}\label{central_weights_2}
\resizebox{14cm}{!}{
\begin{tabular}{|c|c|c|c|c|c|c|c|c|c|c|}
\hline
\textbf{ Alt/M\"{o}bius } & $\mathbf{m(\left\{1\right\})}$ & $\mathbf{m(\left\{2\right\})}$ & $\mathbf{m(\left\{3\right\})}$ & $\mathbf{m(\left\{4\right\})}$ & $\mathbf{m(\left\{1,2\right\})}$ & $\mathbf{m(\left\{1,3\right\})}$ & $\mathbf{m(\left\{1,4\right\})}$ & $\mathbf{m(\left\{2,3\right\})}$ & $\mathbf{m(\left\{2,4\right\})}$ & $\mathbf{m(\left\{3,4\right\})}$ \\
\hline
   $\mathbf{a_{1}}$ &   0.38  & 0.10  & 0.14  & 0.11  & 0.18  & 0.02  & -0.02 & 0.13  & -0.04 & 0.01 \\ 
   $\mathbf{a_{5}}$ & 0.26  & 0.10  & 0.20  & 0.11  & 0.06  & -0.05 & 0.11  & 0.04  & -0.06 & 0.23 \\
   $\mathbf{a_{7}}$ & 0.32  & 0.22  & 0.18  & 0.15  & 0.12  & -0.02 & 0.06  & 0.05  & -0.08 & 0.00 \\
    $\mathbf{a_{11}}$ & 0.35  & 0.10  & 0.20  & 0.16  & 0.07  & 0.08  & 0.00  & 0.06  & -0.05 & 0.03 \\
 %   $a_{12}$ & 0.00  & 0.00  & 0.00  & 0.00  & 0.00  & 0.00  & 0.00  & 0.00  & 0.00  & 0.00 \\
    $\mathbf{a_{13}}$ & 0.44  & 0.07  & 0.16  & 0.08  & 0.23  & -0.08 & -0.01 & 0.05  & -0.03 & 0.08 \\
    $\mathbf{a_{15}}$ & 0.52  & 0.10  & 0.31  & 0.17  & 0.06  & -0.20 & 0.03  & 0.07  & -0.05 & -0.01 \\
    $\mathbf{a_{16}}$ & 0.31  & 0.07  & 0.39  & 0.12  & 0.04  & -0.20 & 0.22  & 0.02  & -0.04 & 0.07 \\
    $\mathbf{a_{17}}$ & 0.24  & 0.10  & 0.11  & 0.10  & 0.16  & 0.10  & 0.06  & 0.13  & -0.04 & 0.04 \\
   \hline
\end{tabular}}  
 
\end{center}
\end{table}

%%%%%%%%%%%%%%%%%%%%%%%%%%%%%%%%%%%%%%%%%%%%%%%%%%%%%%%%%%%%%%%%%%%%%%%%%%%%%%%%
%%%%%%%%%%%%%%%%%%%%%%%%%%%%%%%%%%%%%%%%%%%%%%%%%%%%%%%%%%%%%%%%%%%%%%%%%%%%%%%%

\begin{table}[h]
\begin{center}
\caption{ M\"{o}bius representation of the barycenter of compatible  capacities taking into account precise evaluations on considered criteria and preference information both in terms of comparison of alternatives and importance and interaction of criteria}\label{barycent_2}
\resizebox{14cm}{!}{\begin{tabular}{|c|c|c|c|c|c|c|c|c|c|}
\hline
$\mathbf{m(\left\{1\right\})}$ & $\mathbf{m(\left\{2\right\})}$ & $\mathbf{m(\left\{3\right\})}$ & $\mathbf{m(\left\{4\right\})}$ & $\mathbf{m(\left\{1,2\right\})}$ & $\mathbf{m(\left\{1,3\right\})}$ & $\mathbf{m(\left\{1,4\right\})}$ & $\mathbf{m(\left\{2,3\right\})}$ & $\mathbf{m(\left\{2,4\right\})}$ & $\mathbf{m(\left\{3,4\right\})}$ \\ 
\hline
0.39 & 0.12 & 0.23 & 0.15 & 0.086 & -0.045 & 0.033 & 0.067 & -0.053 & 0.02 \\
\hline  
\end{tabular}}  
\end{center}
\end{table}

%%%%%%%%%%%%%%%%%%%%%%%%%%%%%%%%%%%%%%%%%%%%%%%%%%%%%%%%%%%%%%%%%%%%%%%%%%%%%%%%
%%%%%%%%%%%%%%%%%%%%%%%%%%%%%%%%%%%%%%%%%%%%%%%%%%%%%%%%%%%%%%%%%%%%%%%%%%%%%%%%

\begin{table}[h]
\begin{center}
\caption{Frequency of the preference between pairs of alternatives taking into account precise evaluations on considered criteria  and preference information both in terms of comparison of alternatives and importance and interaction of criteria}\label{preference_2}
\resizebox{17cm}{!}{\begin{tabular}{|c|c|c|c|c|c|c|c|c|c|c|c|c|c|c|c|c|c|c|}
\hline
\textbf{Alt/Alt}   & $\mathbf{a_{1}}$ & $\mathbf{a_{2}}$ & $\mathbf{a_{3}}$ & $\mathbf{a_{4}}$ & $\mathbf{a_{5}}$ & $\mathbf{a_{6}}$ & $\mathbf{a_{7}}$ & $\mathbf{a_{8}}$ & $\mathbf{a_{9}}$ & $\mathbf{a_{10}}$ & $\mathbf{a_{11}}$ & $\mathbf{a_{12}}$ & $\mathbf{a_{13}}$ & $\mathbf{a_{14}}$ & $\mathbf{a_{15}}$ & $\mathbf{a_{16}}$ & $\mathbf{a_{17}}$ & $\mathbf{a_{18}}$ \\
\hline
   $\mathbf{a_{1}}$ & 0.00  & 90.34 & 60.44 & 99.91 & 71.66 & 99.66 & 57.00 & 98.28 & 100.00 & 100.00 & 24.60 & 99.68 & 56.53 & 99.71 & 41.99 & 74.29 & 36.57 & 100.00 \\
   $\mathbf{a_{2}}$ & 9.66  & 0.00  & 19.32 & 79.09 & 2.42  & 94.84 & 15.90 & 46.24 & 99.94 & 96.38 & 1.35  & 83.35 & 12.29 & 88.01 & 12.78 & 0.00  & 6.29  & 60.07 \\
   $\mathbf{a_{3}}$ &39.56 & 80.68 & 0.00  & 97.31 & 58.69 & 99.86 & 45.88 & 84.21 & 100.00 & 98.58 & 28.08 & 98.54 & 36.69 & 100.00 & 2.56  & 62.93 & 37.81 & 87.02 \\
   $\mathbf{a_{4}}$ & 0.09  & 20.91 & 2.69  & 0.00  & 6.07  & 70.65 & 0.05  & 14.80 & 100.00 & 93.23 & 0.09  & 57.33 & 0.42  & 53.37 & 2.23  & 8.81  & 0.03  & 19.30 \\
   $\mathbf{a_{5}}$ & 28.34 & 97.58 & 41.31 & 93.93 & 0.00  & 99.73 & 37.63 & 81.39 & 99.99 & 98.96 & 6.25  & 97.35 & 35.18 & 99.26 & 30.14 & 77.75 & 20.32 & 89.19 \\
   $\mathbf{a_{6}}$ & 0.34  & 5.16  & 0.14  & 29.35 & 0.27  & 0.00  & 1.77  & 10.55 & 97.88 & 73.47 & 0.00  & 22.66 & 0.85  & 6.08  & 0.00  & 0.39  & 0.17  & 14.21 \\
   $\mathbf{a_{7}}$ & 43.00 & 84.10 & 54.12 & 99.95 & 62.37 & 98.23 & 0.00  & 91.99 & 100.00 & 100.00 & 32.00 & 99.10 & 46.98 & 97.86 & 38.11 & 65.33 & 39.23 & 97.40 \\
   $\mathbf{a_{8}}$ & 1.72  & 53.76 & 15.79 & 85.20 & 18.61 & 89.45 & 8.01  & 0.00  & 100.00 & 99.06 & 0.41  & 85.96 & 0.00  & 84.80 & 9.68  & 22.92 & 1.66  & 68.45 \\
   $\mathbf{a_{9}}$ & 0.00  & 0.06  & 0.00  & 0.00  & 0.01  & 2.12  & 0.00  & 0.00  & 0.00  & 0.99  & 0.00  & 0.00  & 0.00  & 0.32  & 0.00  & 0.02  & 0.00  & 0.00 \\
   $\mathbf{a_{10}}$ & 0.00  & 3.62  & 1.42  & 6.77  & 1.04  & 26.53 & 0.00  & 0.94  & 99.01 & 0.00  & 0.00  & 11.80 & 0.02  & 14.84 & 1.04  & 1.89  & 0.00  & 1.00 \\
   $\mathbf{a_{11}}$ & 75.40 & 98.65 & 71.92 & 99.91 & 93.75 & 100.00 & 68.00 & 99.59 & 100.00 & 100.00 & 0.00  & 100.00 & 71.33 & 100.00 & 53.63 & 95.71 & 72.88 & 100.00 \\
   $\mathbf{a_{12}}$ & 0.32  & 16.65 & 1.46  & 42.67 & 2.65  & 77.34 & 0.90  & 14.04 & 100.00 & 88.20 & 0.00  & 0.00  & 0.97  & 50.18 & 1.09  & 4.92  & 0.10  & 18.24 \\
   $\mathbf{a_{13}}$ & 43.47 & 87.71 & 63.31 & 99.58 & 64.82 & 99.15 & 53.02 & 100.00 & 100.00 & 99.98 & 28.67 & 99.03 & 0.00  & 98.96 & 33.03 & 68.69 & 39.70 & 98.89 \\
   $\mathbf{a_{14}}$ & 0.29  & 11.99 & 0.00  & 46.63 & 0.74  & 93.92 & 2.14  & 15.20 & 99.68 & 85.16 & 0.00  & 49.82 & 1.04  & 0.00  & 0.40  & 1.66  & 0.07  & 19.28 \\
   $\mathbf{a_{15}}$ & 58.01 & 87.22 & 97.44 & 97.77 & 69.86 & 100.00 & 61.89 & 90.32 & 100.00 & 98.96 & 46.37 & 98.91 & 66.97 & 99.60 & 0.00  & 73.79 & 53.59 & 90.86 \\
   $\mathbf{a_{16}}$ & 25.71 & 100.00 & 37.07 & 91.19 & 22.25 & 99.61 & 34.67 & 77.08 & 99.98 & 98.11 & 4.29  & 95.08 & 31.31 & 98.34 & 26.21 & 0.00  & 18.58 & 83.79 \\
   $\mathbf{a_{17}}$ & 63.43 & 93.71 & 62.19 & 99.97 & 79.68 & 99.83 & 60.77 & 98.34 & 100.00 & 100.00 & 27.12 & 99.90 & 60.30 & 99.93 & 46.41 & 81.42 & 0.00  & 100.00 \\
   $\mathbf{a_{18}}$ & 0.00  & 39.93 & 12.98 & 80.70 & 10.81 & 85.79 & 2.60  & 31.55 & 100.00 & 99.00 & 0.00  & 81.76 & 1.11  & 80.72 & 9.14  & 16.21 & 0.00  & 0.00 \\
\hline
\end{tabular}} 
 \end{center}
\end{table}

%%%%%%%%%%%%%%%%%%%%%%%%%%%%%%%%%%%%%%%%%%%%%%%%%%%%%%%%%%%%%%% 
%% IMPRECISE EVALUATIONS AND PREFERENCES ONLY ON CRITERIA %%%%%
%%%%%%%%%%%%%%%%%%%%%%%%%%%%%%%%%%%%%%%%%%%%%%%%%%%%%%%%%%%%%%%

\begin{table}[h]
\begin{center}
\caption{Rank acceptability indices  taking into account imprecise evaluations and preference information both in terms of comparison of alternatives and importance and interaction of criteria}\label{rank_accept_IMP_1}
\resizebox{17cm}{!}{
\begin{tabular}{|c|c|c|c|c|c|c|c|c|c|c|c|c|c|c|c|c|c|c|}
\hline
\textbf{Alt} & $\mathbf{b_k^{1}}$ & $\mathbf{b_k^{2}}$ & $\mathbf{b_k^{3}}$ & $\mathbf{b_k^{4}}$ & $\mathbf{b_k^{5}}$ & $\mathbf{b_k^{6}}$ & $\mathbf{b_k^{7}}$ & $\mathbf{b_k^{8}}$ & $\mathbf{b_k^{9}}$ & $\mathbf{b_k^{10}}$ & $\mathbf{b_k^{11}}$ & $\mathbf{b_k^{12}}$ & $\mathbf{b_k^{13}}$ & $\mathbf{b_k^{14}}$ & $\mathbf{b_k^{15}}$ & $\mathbf{b_k^{16}}$ & $\mathbf{b_k^{17}}$ & $\mathbf{b_k^{18}}$ \\  
\hline

   $\mathbf{a_1}$    &   7.80  & 11.66 & 12.23 & 10.72 & 11.48 & 10.95 & 9.26  & 7.45  & 5.23  & 3.89  & 3.05  & 2.22  & 1.63  & 1.25  & 0.83  & 0.33  & 0.03  & 0.00 \\
   $\mathbf{a_2}$    &   9.42  & 7.08  & 7.25  & 7.84  & 8.24  & 8.19  & 7.56  & 7.69  & 6.93  & 7.52  & 6.19  & 5.50  & 3.96  & 2.98  & 2.36  & 1.15  & 0.11  & 0.01 \\
   $\mathbf{a_3}$    &   2.16  & 4.16  & 4.12  & 4.34  & 4.94  & 5.64  & 6.01  & 5.69  & 6.05  & 6.32  & 7.07  & 7.76  & 7.43  & 7.82  & 8.40  & 7.87  & 3.97  & 0.26 \\
   $\mathbf{a_4}$    &  0.10  & 0.29  & 0.51  & 0.82  & 1.20  & 1.55  & 2.07  & 3.25  & 3.99  & 5.40  & 6.84  & 9.38  & 12.79 & 14.08 & 15.01 & 14.78 & 7.96  & 0.00 \\
   $\mathbf{a_5}$    &   6.34  & 9.60  & 10.55 & 9.78  & 8.52  & 7.84  & 7.50  & 8.30  & 7.63  & 6.38  & 5.68  & 4.53  & 3.32  & 2.58  & 1.12  & 0.30  & 0.02  & 0.00 \\
   $\mathbf{a_6}$    & 0.50  & 0.92  & 1.25  & 1.36  & 1.69  & 2.23  & 2.81  & 3.25  & 3.92  & 4.77  & 4.90  & 5.47  & 6.82  & 8.83  & 11.15 & 20.99 & 17.64 & 1.49 \\
   $\mathbf{a_7}$    & 11.43 & 6.15  & 5.87  & 5.83  & 5.74  & 6.12  & 6.89  & 6.70  & 7.95  & 8.05  & 7.38  & 6.97  & 5.79  & 4.14  & 3.07  & 1.71  & 0.20  & 0.01 \\
   $\mathbf{a_8}$    & 1.87  & 2.58  & 3.53  & 4.74  & 6.48  & 8.21  & 9.03  & 9.63  & 9.40  & 9.53  & 9.50  & 6.99  & 5.38  & 4.53  & 4.30  & 3.26  & 0.65  & 0.38 \\
   $\mathbf{a_9}$    &  0.00  & 0.00  & 0.00  & 0.00  & 0.00  & 0.00  & 0.00  & 0.01  & 0.02  & 0.04  & 0.07  & 0.12  & 0.20  & 0.34  & 0.74  & 1.23  & 7.17  & 90.05 \\
   $\mathbf{a_{10}}$ &  0.00  & 0.00  & 0.00  & 0.00  & 0.00  & 0.02  & 0.06  & 0.15  & 0.49  & 1.04  & 1.84  & 2.65  & 4.49  & 6.51  & 8.41  & 16.19 & 51.76 & 6.37 \\
   $\mathbf{a_{11}}$ & 21.06 & 16.32 & 15.09 & 13.48 & 10.60 & 7.99  & 6.35  & 3.62  & 2.28  & 1.38  & 0.82  & 0.50  & 0.28  & 0.15  & 0.05  & 0.01  & 0.00  & 0.00 \\
   $\mathbf{a_{12}}$ & 1.64  & 1.92  & 1.84  & 1.84  & 2.09  & 2.57  & 3.27  & 4.31  & 5.31  & 6.53  & 8.42  & 9.65  & 11.80 & 14.47 & 14.38 & 8.17  & 1.80  & 0.00 \\
   $\mathbf{a_{13}}$ & 3.17  & 5.75  & 6.78  & 7.65  & 6.95  & 6.69  & 6.79  & 6.76  & 7.56  & 7.42  & 6.44  & 5.67  & 5.19  & 5.28  & 4.89  & 4.70  & 1.77  & 0.55 \\
   $\mathbf{a_{14}}$ &  3.34  & 3.08  & 2.99  & 3.28  & 4.01  & 5.03  & 5.79  & 6.23  & 7.87  & 7.69  & 8.33  & 8.50  & 11.74 & 9.55  & 7.79  & 4.14  & 0.62  & 0.00 \\
   $\mathbf{a_{15}}$ & 7.53  & 5.50  & 4.79  & 4.78  & 5.26  & 5.55  & 5.50  & 5.36  & 5.64  & 5.95  & 6.38  & 6.61  & 6.30  & 6.41  & 6.68  & 6.29  & 4.96  & 0.53 \\
   $\mathbf{a_{16}}$ & 7.60  & 9.15  & 9.20  & 8.66  & 7.78  & 7.08  & 6.84  & 7.37  & 7.25  & 6.09  & 5.51  & 4.89  & 4.01  & 3.61  & 3.49  & 1.29  & 0.16  & 0.03 \\
   $\mathbf{a_{17}}$ & 15.98 & 15.49 & 13.02 & 12.88 & 11.65 & 9.29  & 7.02  & 4.66  & 3.43  & 2.49  & 1.65  & 1.05  & 0.71  & 0.44  & 0.20  & 0.05  & 0.01  & 0.00 \\
   $\mathbf{a_{18}}$ & 0.06  & 0.36  & 0.98  & 2.03  & 3.37  & 5.06  & 7.26  & 9.56  & 9.06  & 9.52  & 9.91  & 11.53 & 8.15  & 7.02  & 7.13  & 7.54  & 1.17  & 0.31 \\

\hline
\end{tabular}  }
\end{center}
\end{table}

%%%%%%%%%%%%%%%%%%%%%%%%%%%%%%%%%%%%%%%%%%%%%%%%%%%%%%%%%%%%%%%%
%%%%%%%%%%%%%%%%%%%%%%%%%%%%%%%%%%%%%%%%%%%%%%%%%%%%%%%%%%%%%%%%

\begin{table}[h]
\begin{center}
\caption{M\"{o}bius representation of central capacities taking into account imprecise evaluations of alternatives on considered criteria and including preference information in terms of importance and interaction of criteria}\label{central_weights_IMP_1}
\resizebox{14cm}{!}{
\begin{tabular}{|c|c|c|c|c|c|c|c|c|c|c|}
\hline
\textbf{Alt/M\"{o}bius} & $\mathbf{m(\left\{1\right\})}$ & $\mathbf{m(\left\{2\right\})}$ & $\mathbf{m(\left\{3\right\})}$ & $\mathbf{m(\left\{4\right\})}$ & $\mathbf{m(\left\{1,2\right\})}$ & $\mathbf{m(\left\{1,3\right\})}$ & $\mathbf{m(\left\{1,4\right\})}$ & $\mathbf{m(\left\{2,3\right\})}$ & $\mathbf{m(\left\{2,4\right\})}$ & $\mathbf{m(\left\{3,4\right\})}$ \\
\hline
   $\mathbf{a_{1}}$ &  0.29  & 0.11  & 0.18  & 0.18  & 0.16  & 0.00  & -0.03 & 0.14  & -0.05 & 0.01 \\
   $\mathbf{a_{2}}$ &  0.26  & 0.11  & 0.37  & 0.30  & 0.08  & -0.09 & -0.07 & 0.08  & -0.07 & 0.03 \\
   $\mathbf{a_{3}}$ &  0.49  & 0.12  & 0.37  & 0.24  & 0.07  & -0.25 & 0.01  & 0.07  & -0.06 & -0.06 \\
   $\mathbf{a_{4}}$ &  0.24  & 0.16  & 0.32  & 0.41  & 0.11  & -0.03 & -0.02 & 0.04  & -0.03 & -0.20 \\
   $\mathbf{a_{5}}$ &  0.16  & 0.12  & 0.26  & 0.21  & 0.06  & 0.07  & 0.07  & 0.06  & -0.08 & 0.07 \\
   $\mathbf{a_{6}}$ &  0.21  & 0.13  & 0.44  & 0.58  & 0.04  & 0.05  & -0.09 & 0.04  & -0.11 & -0.28 \\
   $\mathbf{a_{7}}$ &  0.30  & 0.20  & 0.22  & 0.21  & 0.18  & -0.08 & 0.02  & 0.06  & -0.09 & -0.01 \\
   $\mathbf{a_{8}}$ &  0.25  & 0.10  & 0.45  & 0.20  & 0.08  & -0.14 & 0.04  & 0.17  & -0.04 & -0.10 \\
   $\mathbf{a_{11}}$ &  0.26  & 0.10  & 0.26  & 0.22  & 0.07  & 0.10  & -0.02 & 0.08  & -0.05 & -0.02 \\
   $\mathbf{a_{12}}$ &  0.24  & 0.13  & 0.37  & 0.51  & 0.10  & -0.01 & -0.09 & 0.06  & -0.07 & -0.25 \\
   $\mathbf{a_{13}}$ &  0.45  & 0.13  & 0.26  & 0.18  & 0.15  & -0.16 & -0.07 & 0.08  & -0.05 & 0.03 \\
   $\mathbf{a_{14}}$ &  0.21  & 0.12  & 0.35  & 0.49  & 0.08  & 0.05  & -0.09 & 0.08  & -0.09 & -0.21 \\
   $\mathbf{a_{15}}$ &  0.50  & 0.12  & 0.36  & 0.24  & 0.07  & -0.23 & -0.02 & 0.08  & -0.06 & -0.04 \\
   $\mathbf{a_{16}}$ &  0.18  & 0.11  & 0.35  & 0.22  & 0.05  & 0.01  & 0.06  & 0.05  & -0.07 & 0.03 \\
   $\mathbf{a_{17}}$ &  0.20  & 0.10  & 0.18  & 0.17  & 0.13  & 0.08  & 0.02  & 0.15  & -0.04 & 0.00 \\
   $\mathbf{a_{18}}$ &  0.16  & 0.10  & 0.24  & 0.17  & 0.11  & 0.06  & 0.05  & 0.23  & -0.06 & -0.06 \\
\hline
\end{tabular}  
 }
\end{center}
\end{table}

%%%%%%%%%%%%%%%%%%%%%%%%%%%%%%%%%%%%%%%%%%%%%%%%%%%%%%%%%%%%%%%%
%%%%%%%%%%%%%%%%%%%%%%%%%%%%%%%%%%%%%%%%%%%%%%%%%%%%%%%%%%%%%%%%

\begin{table}[h]
\begin{center}
\caption{M\"{o}bius representation of the barycenter of the compatible capacities taking into account interval evaluations and preference information of alternatives on considered criteria in terms of importance and interaction of  criteria}\label{barycent_3}
\resizebox{14cm}{!}{\begin{tabular}{|c|c|c|c|c|c|c|c|c|c|}
\hline
$\mathbf{m(\left\{1\right\})}$ & $\mathbf{m(\left\{2\right\})}$ & $\mathbf{m(\left\{3\right\})}$ & $\mathbf{m(\left\{4\right\})}$ & $\mathbf{m(\left\{1,2\right\})}$ & $\mathbf{m(\left\{1,3\right\})}$ & $\mathbf{m(\left\{1,4\right\})}$ & $\mathbf{m(\left\{2,3\right\})}$ & $\mathbf{m(\left\{2,4\right\})}$ & $\mathbf{m(\left\{3,4\right\})}$ \\ 
\hline
0.27 &   0.12 &   0.27 &   0.23 &   0.10 &   -0.007 &   -0.005 &   0.09 &    -0.06 &   -0.013 \\
\hline  
\end{tabular}}  
\end{center}
\end{table}

%%%%%%%%%%%%%%%%%%%%%%%%%%%%%%%%%%%%%%%%%%%%%%%%%%%%%%%%%%%%%%%%
%%%%%%%%%%%%%%%%%%%%%%%%%%%%%%%%%%%%%%%%%%%%%%%%%%%%%%%%%%%%%%%%

\begin{table}[h]
\begin{center}
\caption{Frequency of the preference between pairs of alternatives taking into account imprecise evaluations of alternatives on considered criteria and preference information expressed in terms of  importance and interaction of criteria only }\label{preference_IMP_1}
\resizebox{17cm}{!}{

\begin{tabular}{|c|c|c|c|c|c|c|c|c|c|c|c|c|c|c|c|c|c|c|}
\hline
\textbf{Alt/Alt}   & $\mathbf{a_{1}}$ & $\mathbf{a_{2}}$ & $\mathbf{a_{3}}$ & $\mathbf{a_{4}}$ & $\mathbf{a_{5}}$ & $\mathbf{a_{6}}$ & $\mathbf{a_{7}}$ & $\mathbf{a_{8}}$ & $\mathbf{a_{9}}$ & $\mathbf{a_{10}}$ & $\mathbf{a_{11}}$ & $\mathbf{a_{12}}$ & $\mathbf{a_{13}}$ & $\mathbf{a_{14}}$ & $\mathbf{a_{15}}$ & $\mathbf{a_{16}}$ & $\mathbf{a_{17}}$ & $\mathbf{a_{18}}$ \\
\hline
    $\mathbf{a_{1}}$ &   0.00  & 58.70 & 78.99 & 91.93 & 54.65 & 86.95 & 67.84 & 74.97 & 99.56 & 99.69 & 33.11 & 84.96 & 77.10 & 77.84 & 72.88 & 56.05 & 34.62 & 88.88 \\
    $\mathbf{a_{2}}$ &   41.30 & 0.00  & 64.82 & 84.80 & 44.84 & 91.95 & 53.41 & 63.46 & 99.91 & 98.53 & 26.97 & 80.98 & 57.54 & 74.49 & 61.01 & 47.78 & 32.65 & 73.23 \\
    $\mathbf{a_{3}}$ &  21.01 & 35.18 & 0.00  & 70.30 & 34.29 & 67.41 & 30.46 & 43.17 & 99.11 & 86.97 & 17.40 & 59.37 & 33.71 & 50.44 & 38.55 & 36.42 & 19.39 & 51.98 \\
    $\mathbf{a_{4}}$ & 8.07  & 15.20 & 29.70 & 0.00  & 14.00 & 53.61 & 12.01 & 21.74 & 100.00 & 82.30 & 4.64  & 34.27 & 19.97 & 25.77 & 27.52 & 16.88 & 5.18  & 29.96 \\
    $\mathbf{a_{5}}$ & 45.35 & 55.16 & 65.71 & 86.00 & 0.00  & 92.65 & 56.97 & 68.91 & 99.88 & 99.28 & 22.34 & 81.66 & 60.05 & 77.52 & 61.78 & 56.91 & 34.60 & 79.78 \\
    $\mathbf{a_{6}}$ & 13.05 & 8.05  & 32.59 & 46.39 & 7.35  & 0.00  & 23.11 & 21.81 & 98.43 & 74.52 & 4.75  & 30.60 & 26.35 & 13.00 & 30.55 & 8.61  & 8.25  & 28.71 \\
    $\mathbf{a_{7}}$ & 32.16 & 46.59 & 69.54 & 87.99 & 43.03 & 76.89 & 0.00  & 58.11 & 99.50 & 97.24 & 27.80 & 75.82 & 59.52 & 63.88 & 62.11 & 44.89 & 27.60 & 68.24 \\
    $\mathbf{a_{8}}$ & 25.03 & 36.54 & 56.83 & 78.26 & 31.09 & 78.19 & 41.89 & 0.00  & 99.10 & 97.76 & 12.07 & 70.11 & 47.26 & 60.70 & 52.63 & 33.71 & 14.22 & 69.80 \\
    $\mathbf{a_{9}}$ & 0.44  & 0.09  & 0.89  & 0.00  & 0.12  & 1.57  & 0.50  & 0.90  & 0.00  & 7.34  & 0.08  & 0.01  & 1.49  & 0.10  & 1.04  & 0.26  & 0.16  & 1.15 \\
    $\mathbf{a_{10}}$ & 0.31  & 1.47  & 13.03 & 17.70 & 0.72  & 25.48 & 2.76  & 2.24  & 92.66 & 0.00  & 0.06  & 13.69 & 6.20  & 8.05  & 12.53 & 1.46  & 0.12  & 3.40 \\
    $\mathbf{a_{11}}$ & 66.89 & 73.03 & 82.60 & 95.36 & 77.66 & 95.25 & 72.20 & 87.93 & 99.92 & 99.94 & 0.00  & 91.39 & 77.31 & 88.28 & 77.35 & 77.99 & 57.60 & 95.87 \\
    $\mathbf{a_{12}}$ & 15.04 & 19.02 & 40.63 & 65.73 & 18.34 & 69.40 & 24.18 & 29.89 & 99.99 & 86.31 & 8.61  & 0.00  & 29.68 & 29.69 & 37.78 & 21.38 & 10.72 & 37.54 \\
    $\mathbf{a_{13}}$ & 22.90 & 42.46 & 66.29 & 80.03 & 39.95 & 73.65 & 40.48 & 52.74 & 98.51 & 93.80 & 22.69 & 70.32 & 0.00  & 59.91 & 56.87 & 42.00 & 22.21 & 62.07 \\
    $\mathbf{a_{14}}$ & 22.16 & 25.51 & 49.56 & 74.23 & 22.48 & 87.00 & 36.12 & 39.30 & 99.90 & 91.95 & 11.72 & 70.31 & 40.09 & 0.00  & 46.50 & 26.17 & 15.16 & 48.14 \\
    $\mathbf{a_{15}}$ & 27.12 & 38.99 & 60.83 & 72.48 & 38.22 & 69.45 & 37.89 & 47.37 & 98.96 & 87.47 & 22.65 & 62.22 & 43.13 & 53.50 & 0.00  & 40.20 & 24.45 & 55.56 \\
    $\mathbf{a_{16}}$ & 43.95 & 52.22 & 63.58 & 83.12 & 43.09 & 91.39 & 55.11 & 66.29 & 99.74 & 98.54 & 22.01 & 78.62 & 58.00 & 73.83 & 59.80 & 0.00  & 34.14 & 75.92 \\
    $\mathbf{a_{17}}$ & 65.29 & 67.35 & 80.61 & 94.82 & 65.40 & 91.75 & 72.40 & 85.78 & 99.84 & 99.88 & 42.40 & 89.28 & 77.79 & 84.84 & 75.55 & 65.86 & 0.00  & 97.38 \\
    $\mathbf{a_{18}}$ & 11.12 & 26.77 & 48.02 & 70.04 & 20.22 & 71.29 & 31.76 & 30.20 & 98.85 & 96.60 & 4.13  & 62.46 & 37.93 & 51.86 & 44.44 & 24.08 & 2.62  & 0.00 \\
\hline
\end{tabular}  
}
\end{center}
\end{table}

%%%%%%%%%%%%%%%%%%%%%%%%%%%%%%%%%%%%%%%%%%%%%%%%%%%%%%%%%%%%%%%%%%%%%%%%%%%
%% IMPRECISE EVALUATIONS AND PREFERENCES ON CRITERIA AND ALTERNATIVES %%%%%
%%%%%%%%%%%%%%%%%%%%%%%%%%%%%%%%%%%%%%%%%%%%%%%%%%%%%%%%%%%%%%%%%%%%%%%%%%%

\begin{table}[h]
\begin{center}
\caption{Rank acceptability indices taking into account imprecise evaluations of alternatives on considered criteria, preference information in terms of importance and interaction of criteria and comparisons between alternatives}\label{rank_accept_IMP_2}
\resizebox{17cm}{!}{
\begin{tabular}{|c|c|c|c|c|c|c|c|c|c|c|c|c|c|c|c|c|c|c|}
\hline
\textbf{Alt} & $\mathbf{b_k^{1}}$ & $\mathbf{b_k^{2}}$ & $\mathbf{b_k^{3}}$ & $\mathbf{b_k^{4}}$ & $\mathbf{b_k^{5}}$ & $\mathbf{b_k^{6}}$ & $\mathbf{b_k^{7}}$ & $\mathbf{b_k^{8}}$ & $\mathbf{b_k^{9}}$ & $\mathbf{b_k^{10}}$ & $\mathbf{b_k^{11}}$ & $\mathbf{b_k^{12}}$ & $\mathbf{b_k^{13}}$ & $\mathbf{b_k^{14}}$ & $\mathbf{b_k^{15}}$ & $\mathbf{b_k^{16}}$ & $\mathbf{b_k^{17}}$ & $\mathbf{b_k^{18}}$ \\  
\hline
$\mathbf{a_1}$ &    19.03 & 22.62 & 23.26 & 15.84 & 9.28  & 4.99  & 2.73  & 1.29  & 0.62  & 0.22  & 0.09  & 0.03  & 0.01  & 0.00  & 0.00  & 0.00  & 0.00  & 0.00 \\
$\mathbf{a_2}$ &    0.00  & 0.01  & 0.04  & 0.13  & 0.43  & 1.60  & 3.84  & 7.12  & 10.46 & 14.90 & 17.97 & 19.02 & 12.43 & 7.38  & 4.10  & 0.55  & 0.02  & 0.00 \\
$\mathbf{a_3}$ &    0.63  & 1.55  & 2.56  & 3.93  & 6.36  & 9.95  & 13.55 & 13.66 & 12.51 & 11.35 & 11.17 & 9.05  & 2.88  & 0.80  & 0.05  & 0.00  & 0.00  & 0.00 \\
$\mathbf{a_4}$ &     0.00  & 0.01  & 0.04  & 0.10  & 0.26  & 0.64  & 1.28  & 2.63  & 3.55  & 4.89  & 6.33  & 8.90  & 14.84 & 20.53 & 23.38 & 10.84 & 1.78  & 0.00 \\
$\mathbf{a_5}$ &     0.47  & 1.62  & 3.70  & 7.17  & 10.58 & 12.19 & 12.14 & 12.65 & 12.70 & 10.04 & 7.57  & 5.02  & 2.57  & 1.13  & 0.41  & 0.03  & 0.00  & 0.00 \\
$\mathbf{a_6}$ &     0.00  & 0.00  & 0.00  & 0.00  & 0.00  & 0.01  & 0.02  & 0.04  & 0.08  & 0.19  & 0.43  & 0.84  & 1.96  & 4.75  & 10.73 & 44.11 & 36.42 & 0.41 \\
$\mathbf{a_7}$ &    22.46 & 15.47 & 15.55 & 18.16 & 11.48 & 6.59  & 4.04  & 2.61  & 1.75  & 0.96  & 0.51  & 0.25  & 0.10  & 0.04  & 0.01  & 0.00  & 0.00  & 0.00 \\
$\mathbf{a_8}$ &    0.00  & 0.01  & 0.09  & 0.52  & 2.14  & 7.74  & 11.46 & 13.98 & 12.38 & 12.80 & 11.37 & 9.61  & 7.85  & 5.49  & 3.53  & 1.02  & 0.01  & 0.00 \\
$\mathbf{a_9}$ &    0.00  & 0.00  & 0.00  & 0.00  & 0.00  & 0.00  & 0.00  & 0.00  & 0.00  & 0.00  & 0.00  & 0.00  & 0.00  & 0.00  & 0.00  & 0.01  & 0.59  & 99.40 \\
$\mathbf{a_{10}}$ &    0.00  & 0.00  & 0.00  & 0.00  & 0.00  & 0.00  & 0.00  & 0.00  & 0.01  & 0.03  & 0.07  & 0.17  & 0.60  & 1.91  & 5.88  & 30.78 & 60.35 & 0.19 \\
$\mathbf{a_{11}}$ &   25.30 & 22.25 & 18.37 & 14.01 & 9.37  & 5.43  & 3.37  & 1.37  & 0.40  & 0.10  & 0.03  & 0.01  & 0.00  & 0.00  & 0.00  & 0.00  & 0.00  & 0.00 \\
$\mathbf{a_{12}}$ &    0.00  & 0.00  & 0.00  & 0.02  & 0.05  & 0.15  & 0.40  & 1.04  & 2.07  & 3.63  & 5.81  & 9.13  & 16.47 & 26.31 & 27.59 & 6.90  & 0.44  & 0.00 \\
$\mathbf{a_{13}}$ &    2.29  & 6.10  & 10.01 & 15.94 & 23.21 & 17.91 & 11.67 & 6.41  & 3.78  & 1.70  & 0.66  & 0.22  & 0.07  & 0.02  & 0.00  & 0.00  & 0.00  & 0.00 \\
$\mathbf{a_{14}}$ &    0.00  & 0.00  & 0.00  & 0.00  & 0.02  & 0.09  & 0.39  & 1.51  & 3.35  & 5.66  & 8.71  & 13.71 & 23.24 & 20.85 & 17.91 & 4.23  & 0.32  & 0.00 \\
$\mathbf{a_{15}}$ &    2.91  & 4.15  & 5.40  & 7.43  & 10.53 & 13.45 & 12.39 & 9.96  & 8.22  & 7.05  & 6.43  & 5.05  & 3.38  & 2.07  & 1.24  & 0.31  & 0.05  & 0.00 \\
$\mathbf{a_{16}}$ &    0.47  & 1.03  & 1.69  & 2.93  & 5.28  & 8.65  & 12.29 & 15.63 & 17.18 & 14.82 & 10.75 & 5.60  & 2.60  & 0.99  & 0.07  & 0.00  & 0.00  & 0.00 \\
$\mathbf{a_{17}}$ &   26.42 & 25.15 & 19.06 & 12.73 & 7.65  & 4.48  & 2.69  & 1.20  & 0.46  & 0.12  & 0.03  & 0.01  & 0.00  & 0.00  & 0.00  & 0.00  & 0.00  & 0.00 \\
$\mathbf{a_{18}}$ &    0.00  & 0.03  & 0.23  & 1.10  & 3.36  & 6.12  & 7.74  & 8.91  & 10.48 & 11.55 & 12.07 & 13.37 & 10.98 & 7.74  & 5.08  & 1.21  & 0.02  & 0.00 \\
\hline
\end{tabular}  }
\end{center}
\end{table}

%%%%%%%%%%%%%%%%%%%%%%%%%%%%%%%%%%%%%%%%%%%%%%%%%%%%%%%%%%%%%%%%
%%%%%%%%%%%%%%%%%%%%%%%%%%%%%%%%%%%%%%%%%%%%%%%%%%%%%%%%%%%%%%%%

\begin{table}[h]
\begin{center}
\caption{M\"{o}bius representation of central capacities for alternatives taking into account imprecise evaluations of alternatives on considered criteria, preferences on importance and interaction of criteria and comparisons between alternatives}\label{central_weights_IMP_2}
\resizebox{14cm}{!}{
\begin{tabular}{|c|c|c|c|c|c|c|c|c|c|c|}
\hline
\textbf{Alt/M\"{o}bius} & $\mathbf{m(\left\{1\right\})}$ & $\mathbf{m(\left\{2\right\})}$ & $\mathbf{m(\left\{3\right\})}$ & $\mathbf{m(\left\{4\right\})}$ & $\mathbf{m(\left\{1,2\right\})}$ & $\mathbf{m(\left\{1,3\right\})}$ & $\mathbf{m(\left\{1,4\right\})}$ & $\mathbf{m(\left\{2,3\right\})}$ & $\mathbf{m(\left\{2,4\right\})}$ & $\mathbf{m(\left\{3,4\right\})}$ \\
\hline
$\mathbf{a_{1}}$ &    0.31  & 0.18  & 0.19  & 0.19  & 0.10  & 0.03  & 0.00  & 0.10  & -0.09 & -0.01 \\
$\mathbf{a_{3}}$ &    0.40  & 0.15  & 0.26  & 0.21  & 0.09  & -0.11 & 0.03  & 0.09  & -0.06 & -0.05 \\
$\mathbf{a_{4}}$ &     0.27  & 0.14  & 0.27  & 0.24  & 0.08  & -0.01 & 0.01  & 0.07  & 0.06  & -0.15 \\
$\mathbf{a_{5}}$ &      0.28  & 0.16  & 0.20  & 0.21  & 0.06  & 0.05  & 0.04  & 0.06  & -0.11 & 0.03 \\
$\mathbf{a_{7}}$ &    0.32  & 0.20  & 0.21  & 0.20  & 0.11  & -0.04 & 0.01  & 0.09  & -0.09 & -0.01 \\
$\mathbf{a_{11}}$ &    0.30  & 0.16  & 0.17  & 0.19  & 0.09  & 0.09  & 0.01  & 0.08  & -0.08 & -0.02 \\
$\mathbf{a_{12}}$ &    0.30  & 0.10  & 0.30  & 0.30  & 0.07  & -0.03 & -0.02 & 0.09  & 0.04  & -0.15 \\
$\mathbf{a_{13}}$ &    0.34  & 0.19  & 0.21  & 0.19  & 0.11  & -0.04 & 0.00  & 0.09  & -0.09 & -0.01 \\
$\mathbf{a_{15}}$ &     0.39  & 0.15  & 0.25  & 0.21  & 0.09  & -0.09 & 0.02  & 0.09  & -0.07 & -0.04 \\
$\mathbf{a_{16}}$ &    0.31  & 0.16  & 0.25  & 0.22  & 0.06  & -0.03 & 0.04  & 0.06  & -0.11 & 0.04 \\
$\mathbf{a_{17}}$ &   0.29  & 0.18  & 0.18  & 0.19  & 0.10  & 0.07  & 0.00  & 0.10  & -0.08 & -0.01 \\
$\mathbf{a_{18}}$ &  0.28  & 0.20  & 0.21  & 0.18  & 0.10  & 0.03  & -0.02 & 0.12  & -0.11 & 0.00 \\
\hline
\end{tabular}  
 }
\end{center}
\end{table}

%%%%%%%%%%%%%%%%%%%%%%%%%%%%%%%%%%%%%%%%%%%%%%%%%%%%%%%%%%%%%%%%
%%%%%%%%%%%%%%%%%%%%%%%%%%%%%%%%%%%%%%%%%%%%%%%%%%%%%%%%%%%%%%%%

\begin{table}[h]
\begin{center}
\caption{ M\"{o}bius representation of the barycenter of the compatible capacities taking into account interval evaluations of alternatives on considered criteria, preference information on importance and interaction of criteria and comparisons between alternatives}\label{barycent_4}
\resizebox{14cm}{!}{\begin{tabular}{|c|c|c|c|c|c|c|c|c|c|}
\hline
$\mathbf{m(\left\{1\right\})}$ & $\mathbf{m(\left\{2\right\})}$ & $\mathbf{m(\left\{3\right\})}$ & $\mathbf{m(\left\{4\right\})}$ & $\mathbf{m(\left\{1,2\right\})}$ & $\mathbf{m(\left\{1,3\right\})}$ & $\mathbf{m(\left\{1,4\right\})}$ & $\mathbf{m(\left\{2,3\right\})}$ & $\mathbf{m(\left\{2,4\right\})}$ & $\mathbf{m(\left\{3,4\right\})}$ \\ 
\hline
 0.31 &   0.18 &   0.19 &   0.19 &   0.097 &   0.034 &   0.008 &   0.091 &   -0.08 &   -0.014 \\
\hline  
\end{tabular}}  
\end{center}
\end{table}

%%%%%%%%%%%%%%%%%%%%%%%%%%%%%%%%%%%%%%%%%%%%%%%%%%%%%%%%%%%%%%%%
%%%%%%%%%%%%%%%%%%%%%%%%%%%%%%%%%%%%%%%%%%%%%%%%%%%%%%%%%%%%%%%%

\begin{table}[h]
\begin{center}
\caption{Frequency of the preference between pairs of alternatives taking into account imprecise evaluations of alternatives on considered criteria, preferences on importance and interaction of criteria and comparisons between alternatives }\label{preference_IMP_2}
\resizebox{17cm}{!}{
\begin{tabular}{|c|c|c|c|c|c|c|c|c|c|c|c|c|c|c|c|c|c|c|}
\hline
\textbf{Alt/Alt}   & $\mathbf{a_{1}}$ & $\mathbf{a_{2}}$ & $\mathbf{a_{3}}$ & $\mathbf{a_{4}}$ & $\mathbf{a_{5}}$ & $\mathbf{a_{6}}$ & $\mathbf{a_{7}}$ & $\mathbf{a_{8}}$ & $\mathbf{a_{9}}$ & $\mathbf{a_{10}}$ & $\mathbf{a_{11}}$ & $\mathbf{a_{12}}$ & $\mathbf{a_{13}}$ & $\mathbf{a_{14}}$ & $\mathbf{a_{15}}$ & $\mathbf{a_{16}}$ & $\mathbf{a_{17}}$ & $\mathbf{a_{18}}$ \\
\hline
  $\mathbf{a_{1}}$ &  0.00  & 99.09 & 92.23 & 99.70 & 91.39 & 99.98 & 55.80 & 98.53 & 100.00 & 100.00 & 47.59 & 99.88 & 80.36 & 99.85 & 85.65 & 94.32 & 44.09 & 99.16 \\
  $\mathbf{a_{2}}$ &  0.91  & 0.00  & 29.52 & 78.07 & 14.31 & 99.20 & 2.96  & 37.69 & 99.99 & 99.33 & 0.25  & 83.54 & 6.58  & 77.66 & 23.08 & 0.00  & 0.39  & 45.69 \\
  $\mathbf{a_{3}}$ &  7.77  & 70.48 & 0.00  & 93.07 & 43.69 & 99.57 & 8.24  & 61.99 & 100.00 & 99.87 & 8.77  & 95.18 & 18.31 & 100.00 & 35.82 & 53.00 & 8.02  & 66.50 \\
  $\mathbf{a_{4}}$ &  0.30  & 21.93 & 6.93  & 0.00  & 9.17  & 86.72 & 0.23  & 16.44 & 100.00 & 95.37 & 0.45  & 51.96 & 0.97  & 42.25 & 6.56  & 11.77 & 0.23  & 21.00 \\
  $\mathbf{a_{5}}$ & 8.61  & 85.69 & 56.31 & 90.83 & 0.00  & 99.99 & 15.78 & 65.84 & 100.00 & 99.85 & 2.72  & 94.92 & 27.39 & 94.32 & 46.07 & 62.74 & 5.17  & 73.95 \\
  $\mathbf{a_{6}}$ &  0.02  & 0.80  & 0.43  & 13.28 & 0.01  & 0.00  & 0.14  & 2.68  & 99.59 & 62.13 & 0.00  & 9.16  & 0.13  & 4.64  & 0.78  & 0.04  & 0.00  & 2.85 \\
  $\mathbf{a_{7}}$ & 44.20 & 97.04 & 91.76 & 99.77 & 84.22 & 99.86 & 0.00  & 96.75 & 100.00 & 100.00 & 43.89 & 99.73 & 76.56 & 99.37 & 84.76 & 90.05 & 40.20 & 96.87 \\
  $\mathbf{a_{8}}$ & 1.47  & 62.31 & 38.01 & 83.56 & 34.16 & 97.32 & 3.25  & 0.00  & 100.00 & 99.93 & 2.64  & 86.11 & 0.00  & 79.90 & 29.48 & 42.28 & 1.66  & 58.17 \\
  $\mathbf{a_{9}}$ & 0.00  & 0.01  & 0.00  & 0.00  & 0.00  & 0.41  & 0.00  & 0.00  & 0.00  & 0.20  & 0.00  & 0.00  & 0.00  & 0.00  & 0.00  & 0.00  & 0.00  & 0.00 \\
  $\mathbf{a_{10}}$ & 0.00  & 0.67  & 0.13  & 4.63  & 0.15  & 37.87 & 0.00  & 0.07  & 99.80 & 0.00  & 0.00  & 4.19  & 0.00  & 3.89  & 0.28  & 0.19  & 0.00  & 0.19 \\
  $\mathbf{a_{11}}$ & 52.41 & 99.75 & 91.23 & 99.55 & 97.28 & 100.00 & 56.11 & 97.36 & 100.00 & 100.00 & 0.00  & 99.94 & 77.93 & 99.96 & 84.83 & 97.71 & 47.04 & 98.97 \\
  $\mathbf{a_{12}}$ & 0.12  & 16.46 & 4.82  & 48.04 & 5.08  & 90.84 & 0.27  & 13.89 & 100.00 & 95.81 & 0.06  & 0.00  & 0.60  & 39.11 & 5.38  & 7.14  & 0.04  & 17.50 \\
  $\mathbf{a_{13}}$ &19.64 & 93.42 & 81.69 & 99.03 & 72.61 & 99.87 & 23.44 & 100.00 & 100.00 & 100.00 & 22.07 & 99.40 & 0.00  & 99.07 & 69.73 & 82.77 & 18.51 & 91.79 \\
  $\mathbf{a_{14}}$ &0.15  & 22.34 & 0.00  & 57.75 & 5.68  & 95.36 & 0.63  & 20.10 & 100.00 & 96.11 & 0.04  & 60.89 & 0.93  & 0.00  & 7.73  & 9.83  & 0.05  & 24.04 \\
  $\mathbf{a_{15}}$ & 14.35 & 76.92 & 63.57 & 93.44 & 53.93 & 99.22 & 15.24 & 70.52 & 100.00 & 99.72 & 15.17 & 94.62 & 30.27 & 92.27 & 0.00  & 63.15 & 14.14 & 72.99 \\
  $\mathbf{a_{16}}$ & 5.68  & 100.00 & 47.00 & 88.23 & 37.26 & 99.96 & 9.95  & 57.72 & 100.00 & 99.81 & 2.29  & 92.86 & 17.23 & 90.17 & 36.85 & 0.00  & 4.06  & 65.09 \\
  $\mathbf{a_{17}}$ & 55.82 & 99.61 & 91.98 & 99.77 & 94.83 & 100.00 & 59.80 & 98.34 & 100.00 & 100.00 & 52.96 & 99.96 & 81.49 & 99.95 & 85.86 & 95.94 & 0.00  & 99.46 \\
  $\mathbf{a_{18}}$ & 0.84  & 54.31 & 33.50 & 79.00 & 26.05 & 97.15 & 3.13  & 41.83 & 100.00 & 99.81 & 1.03  & 82.50 & 8.21  & 75.96 & 27.01 & 34.91 & 0.54  & 0.00 \\
\hline
\end{tabular}  
}
\end{center}
\end{table}

%\end{small}

\end{document}